\documentstyle[11pt,amscd,verbatim]{amsart}

\newtheorem{thm}{Theorem}[section]
\newtheorem{prop}[thm]{Proposition}
\newtheorem{lem}[thm]{Lemma}
\newtheorem{cor}[thm]{Corollary}

\theoremstyle{definition}
\newtheorem{defn}[thm]{Definition}
\newtheorem{remark}[thm]{Remark}
\newtheorem{example}[thm]{Example}

\newcommand{\ml}{{\operatorname{m}}}

\newcommand{\bbA}{{\Bbb A}}
\newcommand{\bbG}{{\Bbb G}}

\newcommand{\bbZ}{{\Bbb Z}}

\newcommand{\bbP}{{\Bbb P}}

\newcommand{\cI}{{\cal I}}

\newcommand{\cX}{{\cal X}}

\newcommand{\kbar}{\overline k}

\newcommand{\isomo}{\overset{\sim}{=}}

\newcommand{\brokrarr}{\dasharrow}
\newcommand{\lf}{\mathopen}

\let\r=\mathclose

\newcommand{\notdiv}{\mathrel{\not|}}

\newcommand{\paru}{Parusi\'nski}

\let\to=\longrightarrow

\let\tilde=\widetilde

\numberwithin{equation}{section}

\newcommand{\Ker}{\operatorname{Ker}}

\newcommand{\sdp}{\mathbin{{>}\!{\triangleleft}}}

\newcommand{\GL}{\operatorname{GL}}
\newcommand{\SL}{\operatorname{SL}}
\newcommand{\Spin}{\operatorname{Spin}}
\newcommand{\PGL}{\operatorname{PGL}}
\newcommand{\PGLn}{\PGL_n}
\newcommand{\rank}{\operatorname{rank}}
\newcommand{\Stab}{\operatorname{Stab}}

\newcommand{\Char}{\operatorname{char}} %% \char is already a command
\newcommand{\ed}{\operatorname{ed}}
\newcommand{\Gal}{\operatorname{Gal}}
\newcommand{\galois}{\Gal}

\newcommand{\lra}{\longrightarrow}
\newcommand{\Mat}{\operatorname{M}}
\newcommand{\Mn}{\Mat_n}

\newcommand{\trdeg}{\operatorname{trdeg}}
\newcommand{\Sym}{{\operatorname{S}}}

\newcommand{\Oct}{\operatorname{O}}
\newcommand{\trace}{\tr}

\newsymbol\subsetneq 2328
\newsymbol\onto 1310
\newsymbol\twoheadrightarrow 1310

\begin{document}

\title[G-varieties, 8-30-99]{Essential dimensions of algebraic groups
and a resolution theorem for $G$-varieties}
\author[Z. REICHSTEIN and B. YOUSSIN,  8-28-99]
{Zinovy Reichstein and Boris Youssin, \\ 
(WITH AN APPENDIX BY J{\'A}NOS KOLL{\'A}R AND ENDRE SZAB{\'O})}
% \address{Department of Mathematics, Oregon State University,
% Corvallis, OR 97331}
\thanks{Z. Reichstein was partially supported by NSF grant DMS-9801675}
% \email{zinovy@@math.orst.edu}
% \address{Department of Mathematics and Computer Science,
% University of the Negev, Be'er Sheva', Israel}
% \email{youssin@@math.bgu.ac.il}

% \begin{center}
% WITH AN APPENDIX BY J{\'A}NOS KOLL{\'A}R AND ENDRE SZAB{\'O}
% \end{center}

\begin{abstract}
Let $G$ be an algebraic group and let $X$ be a generically free $G$-variety.
We show that $X$ can be transformed,
by a sequence of blowups with smooth $G$-equivariant centers, into 
a $G$-variety $X'$ with the following property: the stabilizer
of every point of $X'$ is isomorphic to a semidirect product 
$U \sdp A$ of a unipotent group $U$ 
and a diagonalizable group $A$.

As an application of this result, we prove new lower 
bounds on essential dimensions of some algebraic groups. 
We also show that certain polynomials in one variable
cannot be simplified by a Tschirnhaus transformation.
\end{abstract}

\maketitle
\tableofcontents

\section{Introduction}

Let $k$ be an algebraically closed base field of characteristic zero, 
let $G$ be an algebraic group and let $X$ be 
a $G$-variety, both defined over $k$. Assume $X$ is generically 
free, i.e., the $G$-action is free on a dense open subset of $X$. Recall that
by a theorem of Rosenlicht~\cite{rosenlicht1},~\cite{rosenlicht2} the rational
quotient map $X \brokrarr B$ separates orbits of $X$ in general position; 
in other words, we can think of $X$ as a $G$-torsor over $B$.

We shall say that $X$ is defined in dimension $d$ if there exists a dominant
rational map $X \brokrarr X_1$ of generically free $G$-varieties 
\begin{equation} \label{e.compr}
 \begin{array}{ccc}  X      & \brokrarr & X_1  \\
                     |      &           &  |           \, \,     \\ 
     \llap{$\pi$\;}  |      &           &  | \rlap{\;$\pi_1$}  \, \, \\ 
                \downarrow  &           &  \downarrow   \, \,  \\ 
                    B       & \brokrarr &   B_1 
                                                    \end{array} %\; , 
\end{equation}
with $\dim(B_1)  \leq d$. (Here the vertical arrows represent 
rational quotient maps for the $G$-action.) The smallest
integer $d$ such that $X$ is defined in dimension $d$ will 
be called the {\em essential dimension} of $X$ and denoted 
by $\ed(X)$; cf. Definition~\ref{def.ed}.
In the sequel we shall refer to the rational map \eqref{e.compr}
as a {\em compression} (or a $G$-{\em compression}) of $X$; see~\S\ref{prel2}.

We will say that the essential dimension $\ed(G)$ of the group $G$
is equal to $d$ if every generically free $G$-variety 
is defined in dimension $d$,
and $d$ is the smallest integer with this property.
The essential dimension is a numerical invariant of the group; it  
can often be characterized as the minimal number of independent 
parameters required to describe all algebraic objects
of a certain type. These objects are field extensions 
if $G = \Sym_n$, division algebras if $G = \PGL_n$, quadratic
forms if $G= O_n$, Cayley algebras if $G = G_2$, Albert algebras 
if $G=F_4$, etc. Groups of essential dimension 0
are precisely the {\em special groups} introduced by 
Serre~\cite{serre1} and classified by Grothendieck~\cite{grothendieck}
in the 1950s.  For details we refer the reader to \cite{r2}; for results 
on essential dimensions of finite groups see also \cite{br1} and \cite{br2}.

The lower bounds on $\ed(G)$ in \cite{r2} are proved in one of
two ways. One approach, due to J.-P. Serre, uses cohomological
invariants (see Lemma~\ref{lem.ci} and \cite[Section~12]{r2}); 
the second method, due to 
the first author, relies on applying the Tsen---Lang theorem to
appropriately defined anisotropic forms. 

In this paper we develop an alternative approach, based on the following 
resolution procedure.

\begin{thm} \label{thm1.1} (Corollary~\ref{cor.b1-1} and Theorem~\ref{thm3.2})
Let $X$ be a generically free variety. Then
there exists a sequence
\[
X_n \stackrel{\pi_n}{\lra} X_{n-1}
\dots \stackrel{\pi_2}{\lra} X_1 \stackrel{\pi_1}{\lra} X_0 = X
\]
of blowups with smooth $G$-invariant centers such that $X_n$ is smooth 
and for every $x \in X_n$ the stabilizer $\Stab(x)$ is isomorphic 
to a semidirect 
product $U \sdp A$, where $U$ is unipotent and $A$ is diagonalizable.
\qed
\end{thm}

In fact, we show that a sequence of equivariant blowups can be chosen
so that $X_n$ is in ``standard form"; 
see Definition~\ref{def3.1} and Corollary~\ref{cor.b1-1}. 
The proof of this result
depends on canonical resolution of singularities; see Section~\ref{sect3}. 

In Sections \ref{sect5}--\ref{sect7} we use the above resolution procedure
to prove the following lower bound on $\ed(X)$ and $\ed(G)$, and the related
numerical invariants $\ed(X;p)$ and $\ed(G;p)$; see Definition~\ref{def.edp}. 

\begin{thm} \label{thm1.2ab} 
Let $G$ be a semisimple group and let $H$ be an 
abelian subgroup of $G$, whose centralizer is finite.

\smallskip
(a) (Theorem~\ref{thm5-5}) Suppose
$X$ is a generically free $G$-variety, $x$ is a smooth point of 
$X$, and $\Stab(x)$ contains $H$. Then $\ed(X) \geq \rank(H)$. 
If $H$ is a $p$-group then $\ed(X; p) \geq \rank(H)$. 

\smallskip
(b) (Theorem~\ref{thm5-6}) $\ed(G) \geq \rank(H)$. If $H$ is a $p$-group
then $\ed(G; p) \geq \rank(H)$.
\end{thm}

Informally speaking, under the assumptions of the theorem,
$x$ is an obstruction to compressing $X$ 
(as in \eqref{e.compr}).  Note that while the essential dimension
is a property of $X$ at the generic point, this obstruction depends 
on the presence of 
special geometric points (namely smooth fixed points of $H$).
This explains our use of biregular methods, such as resolution 
of singularities, in what is apriori a birational setting.

In Section~\ref{sect8} we apply Theorem~\ref{thm1.2ab} to a number of specific
groups $G$. The new bounds we obtain are summarized in the following theorem.
Note that $\ed(G) \geq \ed(G; p)$ for any prime $p$; 
see Definition~\ref{def.edp}. 

\begin{thm} \label{thm1.2cd}
\begin{enumerate}
\item (Theorem~\ref{thm6.1-1a}) $\ed(PO_n; 2) \geq n-1$,

\item (Theorem~\ref{thm7.11}) If $n \equiv 0$ or $\pm 1 \pmod 8$ then
$\ed(\Spin_n; 2) \geq [\dfrac{n}{2}]+1$.

\item
(Theorem~\ref{thm6.1-1}(5--6)) 
$\ed(2E_7; 2) \geq 7$, $\ed(E_7; 2) \geq 8$.
Here $2E_7$ and $E_7$ denote, respectively, the simply connected 
and the adjoint groups of type $E_7$.

\item (Theorem~\ref{thm6.1-1}(7--8)) 
$\ed(E_8; 2) \geq 9$, $\ed(E_8; 3) \geq 5$.
\end{enumerate}
\end{thm} 

We remark that the bound of part (2) is known to be sharp for 
$n = 7$, $8$ and $9$ (see~\cite{rostspin} and Remark~\ref{rem.rost})
and that $\ed(2E_7; 2) \leq \ed(2E_7) \leq 9$ (see~\cite{kordonsky}
and Remark~\ref{rem6.1-1}).  Further results on essential
dimensions of specific groups can be found in Section~\ref{sect8}.

Most previously known lower bounds on $\ed(G)$ can be derived from
the existence of cohomological invariants; see Lemma~\ref{lem.ci}
and~\cite[Section 12]{r2}.  The bounds of Theorem~\ref{thm1.2cd} cannot 
be proved in this way at the moment, since the necessary cohomological 
invariants are not known to exist. However, one can view these bounds
(as well as the bound of Theorem~\ref{thm.pgln}) as an indication of 
what cohomological invariants may exist; see Remark~\ref{rem.ci2}. 

\smallskip
In the last section we give an application of Theorem~\ref{thm1.2ab}(a)
to the problem of simplifying polynomials by Tschirnhaus transformations.
Let $F$ be a field and let \[ \alpha(x) =
x^n + a_1 x^{n-1} + \ldots + a_{n-1} x + a_n \] be an
irreducible polynomial over $F$. Recall that a Tschirnhaus transformation
(without auxiliary radicals) is an isomorphism of fields $F[x]/(\alpha(x)) 
\simeq F[t]/(\beta(t))$, where $\beta(t) \in F[t]$ is another 
irreducible monic polynomial of degree $n$. We shall say 
that $\beta(t)$ is obtained from $\alpha(x)$ via this Tschirnhaus 
transformation. In other words, $\beta(t)$ can be obtained
from $\alpha(x)$ in this way if $\beta(t)$ is the minimal polynomial
of a generator of the field extension $F \subset F[x]/(\alpha(x))$. (Note
that all fields in this paper are assumed to contain a copy of the base field
$k$ and all field extensions and isomorphisms are defined over $k$; see
\S\ref{pre0.5}.)

It is shown in 
\cite{br1} that if $a_1, \ldots, a_n$ are algebraically independent over 
$k$, i.e., $\alpha(x)$ is the general polynomial of degree $n$, then
at least $[n/2]$ coefficients of $\beta(t)$ are again algebraically 
independent over $k$. Our main result here is as follows.

\begin{thm} \label{thm1.3} (Theorem~\ref{thm7.1})
Suppose $\dfrac{n}{2} \leq m \leq n-1$, where $m$ and $n$ are positive
integers. 
Let $a_m, \dots, a_n$ be algebraically independent variables over $k$,
$F = k(a_m, \ldots, a_n)$ and $E = F[x]/f(x)$, where 
$$f(x) = x^n + a_m x^{n-m} + \dots + a_{n-1} x + a_n \; . $$
Then any polynomial obtained from $f(x)$ by a Tschirnhaus 
transformation has at least $n-m$ algebraically independent (over $k$)
coefficients.
\end{thm}

Note that $f(x)$ has $n-m+1$ algebraically independent coefficients. However, 
the form with $n-m$ independent coefficients is easily attained by the
substitution $x = \dfrac{\strut a_n}{a_{n-1}}y$; see the proof of 
Theorem~\ref{thm7.1}.  Thus the lower bound of the theorem is, 
indeed, the best possible. 

Throughout this paper we shall work over a base field $k$ of characteristic 
zero. This assumption will be needed when we appeal to
equivariant resolution of singularities, the Levi decomposition 
of an algebraic group, and the Luna slice theorem. 
We do not know whether or not the results 
of this paper remain valid in prime characteristic.

Theorem~\ref{thm1.1} can be used in various other settings, not directly 
related to compressions or essential dimensions.  In~\cite{ry} we apply it, 
along with the the results of Section~\ref{sect5} and Appendix, to 
the study of splitting fields and splitting groups of $G$-varieties,
including a new construction of noncrossed product division algebras.
In \cite{pl} we apply it give a new algebro-geometric proof 
of the ``Key Lemma'' of \paru~\cite{P}. (The latter result was
used in \paru's proof of the existence of Lipschitz stratifications 
of semianalytic sets.)

We remark that our resolution theorems in Section~\ref{sect3} 
are stated in greater generality than we need for the applications 
given in this paper. In particular, for the sake of these applications, 
it would have sufficed to assume that $k$ is an algebraically closed 
field throughout. (Note, however, that this would not have changed 
the proofs.) The more general statements will be needed for further
applications.

\section*{Acknowledgements}

We would like to thank J.-P. Serre for his help and encouragement.  His 
suggestion to investigate the relationship between the essential dimension 
and the non-toral abelian subgroups of a given algebraic group $G$, was the
starting point for the results of Sections~\ref{sect7}
and~\ref{sect8}.  Serre also contributed Definition~\ref{def.edp}, 
Lemma~\ref{lem.ci}, Remark~\ref{rem.ed=r}, 
the statement of Proposition~\ref{prop6.3} 
and, most importantly, both the statement and the proof 
of Lemma~\ref{lem6.11}.  The last result greatly simplified 
our Theorem~\ref{thm5-6} and subsequent applications. 

We are grateful to
P. D. Milman for many helpful discussions of resolution of singularities,
M. Rost for sharing with us his insights into cohomological invariants,
spin groups and quadratic forms
of low degree, and G. Seitz for answering our questions 
about elementary abelian subgroups of exceptional algebraic groups.

We also thank E. Bierstone, P. D. Milman, M. Rost, and J.-P.
Serre for their comments on earlier versions of this paper.

\section{Notation and terminology}
\label{sect2}

The following notational conventions will be used throughout the paper.
\[
 \begin{array}{lcl}
 k          & & \mbox{a base field of characteristic 0} \\
 \kbar      & & \mbox{the algebraic closure of $k$} \\
 G          & & \mbox{an algebraic group defined over $k$; see~\S\ref{prel1.5}} \\
C(H)= C_G(H)     & & \mbox{the centralizer of $H$ in $G$} \\
\bbA^n=\bbA^n_k & & \mbox{the affine space of dimension $n$ over $k$} \\
\bbG_\ml= \GL_1(k) & & \mbox{the multiplicative group 
$\bbA^1-\{0\}$ over $k$} \\
 X          & & \mbox{an algebraic variety over $k$, often a $G$-variety} \\
 \Stab(x)   & & \mbox{the stabilizer of $x$} \\
\ed    & & \mbox{essential dimension;
                see Definitions~\ref{def.ed} and~\ref{def.edp}} 
\end{array}
\]

\subsection{The base field.} \label{pre0.5} All algebraic objects in this 
paper, such as rings, fields, algebraic groups, algebraic varieties, group
actions, etc. and all maps between them will be defined over a fixed 
base field $k$ of characteristic 0. In Sections~\ref{sect4}--\ref{sect8} we
will generally assume that $k$ is algebraically closed; 
we shall indicate which of the results are true without this
assumption.  In Sections~\ref{sect3} and~\ref{sect9} we will 
not assume that $k$ is algebraically closed.

\subsection{Algebraic varieties.} \label{prel0.5}
Algebraic varieties in this paper
are allowed to be reducible; in other words, an algebraic
variety is a reduced separated scheme of finite type over $k$.  (Note that 
here our terminology is different from that of Hartshorne~\cite{Hart}, 
who defines abstract algebraic varieties to be irreducible.)

Given an algebraic variety $X$, we will denote its ring of rational functions
by $k(X)$, where a rational function on a reducible variety is 
a collection of rational functions on its irreducible components; 
cf.~\S\ref{prel1} below.
Note that $k(X)$ is a field if $X$ is irreducible. In general,
if $X$ has irreducible components $X_i$ then
$k(X)$ is a direct sum of their function fields $k(X_i)$.

Unless otherwise specified, by a point of $X$ we shall always mean 
a closed point.

\subsection{Rational maps} \label{prel1}
A rational map $f \colon X\brokrarr Y$ is an equivalence class of regular
morphisms from dense open subsets of $X$ to $Y$, as in 
\cite[D\'efinition~7.1.2]{ega1}.
Equivalently, $f$ is a collection of rational maps
$f_i\colon X_i\brokrarr Y$, one for each irreducible component $X_i$ of $X$. 
The largest open subset $U$ of $X$ where $f$ is defined is called {\em the 
domain} of $f$; $f(U)$ is called {\em the range} of $f$. A rational map
is said to be {\em dominant} if its range is dense in $Y$.

A dominant rational map $f \colon X\brokrarr Y$  
is said to be $d:1$ if there exists a dense open subset $Y_0$ of its range
such that $f$ is defined on $f^{-1}(Y_0)$ and
$|f^{-1}(y)(\overline{k})| = d$ for every $y \in Y_0(\overline{k})$.

A birational isomorphism between $X$ and $Y$ is a pair of rational maps
$X\brokrarr Y$ and $Y\brokrarr X$ inverse to each other, or equivalently, 
a 1---1 correspondence between the irreducible components $X_i$ of $X$ 
and $Y_i$ of $Y$
and a birational isomorphism between $X_i$ and $Y_i$ for each $i$.

\subsection{Algebraic groups} \label{prel1.5} 
If $G$ is an algebraic group (defined over $k$; see \S\ref{pre0.5})
{\em we shall always assume that $G(k)$
is Zariski dense in $G$}. Note that this is a rather mild assumption; 
in particular, it is obviously satisfied if $k$ is algebraically closed
or if $G$ is a finite group all of whose points are defined over $k$
(e.g., $S_n$, viewed as an algebraic group over $k$).
It is also satisfied if $G$ is connected 
(see~\cite[Theorem~34.4(d)]{humphreys}) and, more generally, if every
irreducible component of G has a $k$-point. 

Our results are, in fact, true, without the above assumption; however, 
leaving it out would complicate the proofs in Section~\ref{sect3} 
(see Remark~\ref{rem.thm.b1-2}). Since this assumption is satisfied 
in every setting we want to consider, we chose to impose it throughout
this paper.

\subsection{$G$-varieties.} \label{prel2}

Let $G$ be an algebraic group.
We shall call an algebraic variety $X$ a $G$-variety if
$X$ is equipped with a regular action of $G$, i.e., an action given by 
a regular morphism $G \times X \lra X$. 

If $X$ and $Y$ are $G$-varieties
then by a regular map $X \lra Y$ of $G$-varieties we mean a regular
$G$-equivariant map. The same applies to rational maps of $G$-varieties,
biregular and birational isomorphisms of $G$-varieties, etc.

A $G$-variety is $X$ called {\em generically free} if $G$ acts freely (i.e.,
with trivial stabilizers) on a dense open subset if $X$.

A {\em $G$-compression} $X \brokrarr Y$ is a dominant rational map of
generically free $G$-varieties. We will also use the term {\em compression}
if the reference to $G$ is clear from the context.

\subsection{Rational quotients and primitive varieties.} \label{prel4}
Let $X$ be a $G$-variety.
A rational map $\pi\colon X \brokrarr Y$ is called {\em the rational quotient 
map} (and $Y$, {\em the rational quotient}) if $\pi^{\ast} (k(Y)) = k(X)^G$.
The rational quotient exists for any $G$-variety; we 
will also denote it by $X/G$.

We will say that $X$ is a {\em primitive} $G$-variety if the rational
quotient $X/G$ is irreducible or, equivalently, if $k(X)^G$ is a field. 
It is easy to see that $X$ is primitive if and only if $G$ transitively
permutes the irreducible components of $X$; see, e.g.,~\cite[Lemma 2.2]{r2}.

By a theorem of Rosenlicht the rational quotient map separates
the $G$-orbits in a dense Zariski open subset of $X$; see 
\cite[Theorem 2]{rosenlicht1}, \cite[Theorem 2.3]{pv} and \cite{rosenlicht2}.
In particular, if $X$ is primitive then each component of $X$ has dimension 
$\dim(Y) + \dim(G)$.
%%%%%%%%%%%%%%%%%%%%%%%%%%%

\section{Equivariant resolution of singularities}
\label{sect3}

Much of this paper relies on the resolution
of singularities theorem and especially on its canonical version which 
only recently became available; see the references below.
In this section we
derive several consequences of this result in the setting of G-varieties.

\begin{defn} \label{def3.1}
We shall say that a generically free $G$-variety $X$ is {\em in standard form
with respect to a divisor} $Y$ if

\smallskip
(i) $X$ is smooth and $Y$ is a normal crossing divisor on $X$

\smallskip
(ii) the $G$-action on $X - Y$ is free, and 

\smallskip
(iii) for every $g \in G$ and for every irreducible component $Z$ of $Y$
either $g(Z) = Z$ or $g(Z) \cap Z = \emptyset$.

\smallskip
We will say that $X$ is {\em in standard form} if it is in standard form
with respect to some divisor $Y$.
\end{defn}

Our interest in $G$-varieties in standard form is explained by the fact
that they have ``small'' stabilizers. This property will be explored 
in Section~\ref{sect4}; see Theorem~\ref{thm3.2}. We will now prove that
every generically free $G$-variety can be brought into standard form
by a sequence of blowups with smooth $G$-equivariant centers.

\begin{thm} \label{thm.b1-2}
Let $X$ be a smooth $G$-variety
and $Y\subset X$ be a closed nowhere dense $G$-invariant
subvariety such that the action of $G$ on $X-Y$ is free. Then
there is a sequence of blowups
\begin{equation} \label{tower1}
\pi\colon X_n \stackrel{\pi_n}{\lra} X_{n-1}
\dots \stackrel{\pi_2}{\lra} X_1 \stackrel{\pi_1}{\lra} X_0 = X 
\end{equation}
with smooth $G$-invariant centers $C_i\subset X_{i}$ such that 
$X_n$ is in standard form with respect to $D_n\cup\pi^{-1}(Y)$, where $D_n$ is
the exceptional divisor of $\pi$ (and, in particular,
$D_n\cup\pi^{-1}(Y)$ is a normal crossing divisor in $X_n$).
\end{thm}

\begin{remark}  \label{rem.thm.b1-2}
Recall that throughout this paper we assume $G(k)$ is Zariski 
dense in $G$; see \S\ref{prel1.5}. This assumption is used only
in this section (in Theorem~\ref{thm.b1-2} and
Corollary~\ref{cor.b1-1})
and only for the purpose of lifting a $G$-action 
on an algebraic variety to its canonical resolution of singularities.
 
In fact, our results
are true without this assumption because {\em an algebraic group action
always lifts to the canonical resolution of singularities}
of Bierstone---Milman~\cite{bm} (see also~\cite{bm1}).

The last assertion follows from the fact that the canonical 
resolution commutes with base field extensions. This
reduces the question of lifting a group action to the case
where $k$ is algebraically closed and thus $G(k)$ is Zariski dense in $G$.
Commutativity with base extensions follows from~\cite[Remark~3.8]{bm}.

Alternatively, the above assertion about lifting the action of $G$
can be derived (by an argument more natural than the one
we give in the proof of Theorem~\ref{thm.b1-2} below)
from the fact that the canonical resolution is functorial 
with respect to smooth morphisms. Functoriality
with respect to smooth morphisms
follows from~\cite[Remark~1.5]{bm} and the constructive
definition of the invariant in~\cite[\S\S4, 6]{bm}.

As we do not need the stronger statements of the results of this section
(without the assumption that $G(k)$ is Zariski dense in $G$), we omit
the details of these arguments.

Note also that it is quite possible that the canonical resolution of
Villamayor~\cite{Vil2} (see also~\cite{Vil}) has the same properties.
\end{remark}

We begin with a preliminary lemma.
Let
\begin{equation} \label{tower2}
\pi\colon X_n \stackrel{\pi_n}{\lra} X_{n-1}
\dots \stackrel{\pi_2}{\lra} X_1 \stackrel{\pi_1}{\lra} X_0 = X 
\end{equation}
be a sequence of blowups with smooth $G$-invariant centers.
Recall that the exceptional divisor $E$ of $\pi$ is the
union of the preimages in $X_n$ of the centers of the blowups
$\pi_1,\dots,\pi_n$; the composition $\pi$ is an isomorphism in the complement
of $E$.
 
\begin{lem} \label{lem4.1-1}
Let $X$ be a $G$-variety, let $\pi\colon X_n \lra X$ be as in \eqref{tower2},
and let $E_1$ be an irreducible component of the exceptional divisor $E$ of
$\pi$.
Then for any $g \in G$, either $g(E_1) = E_1$ or $g(E_1) \cap E_1 = \emptyset$.
\end{lem}

\begin{pf}
Each irreducible component of $E$ is the preimage in $X_n$ of an irreducible
component, say, $C_{i,1}$, of the center $C_i$ of one of the blowups
$\pi_{i+1}\colon X_{i+1}\to X_{i}$.

Since $C_i$ is a smooth $G$-invariant subvariety in $X_{i}$,
its irreducible components $C_{i,1},\dots,C_{i,m}$ are disjoint.

We have $E_1=(\pi_i\dots\pi_n)^{-1}C_{i,1}$; hence,
for any $g\in G$, \[ g(E_1)=(\pi_i\dots\pi_n)^{-1}g(C_{i,1}) \; . \]
As $C_i$ is $G$-invariant and $C_{i,1}$ is its connected component, 
$g(C_{i,1})$
is also a connected component of $C_i$, say, $g(C_{i,1})=C_{i,j}$. Thus 
\[ g(E_1)=(\pi_i\dots\pi_n)^{-1}C_{i,j} \; . \]
If $j=1$ then $g(E_1) = E_1$; if $j\ne1$ then $g(E_1) \cap E_1 = \emptyset$,
since $C_{i,1}$ and $C_{i,j}$ are disjoint.
\end{pf}

\begin{pf*}{Proof of Theorem~\ref{thm.b1-2}}
Let $D_i$ be the exceptional divisor of $\pi_1\dots\pi_{i}\colon X_i\to X$.
Inductively, assume that  $D_{i}$ is a normal crossing divisor in
$X_{i}$. We shall give a construction of each blowup center $C_{i}$ so
that $C_i$ and $D_{i}$ simultaneously have only normal crossings.
It was observed by Hironaka~\cite{hironaka} that this implies that $D_{i+1}$ is
a normal crossing divisor in $X_{i+1}$; this way all $D_i$ are normal crossing
divisors.

Denote by $Y_i$ the union of $D_i$ and the preimage of $Y$ in $X_i$.
The algorithm to choose the blowup centers is as follows.  Let
\begin{equation} \label{eqn5mar1}
X_{l-1}@>\pi_{l-1}>>\dots @>\pi_1>>X_0=X
\end{equation}
be a canonical embedded resolution of singularities of $Y\subset X$, as 
in~\cite[Theorem~1.6]{bm}; then $D_{l-1}$ and the strict transform
$C_{l-1}$ of $Y$ in $X_{l-1}$ simultaneously have only normal crossings.

Let 
\begin{equation} \label{eqn13mar1}
X_l@>\pi_l>>X_{l-1}
\end{equation}
be the blowup centered at $C_{l-1}$; then $Y_l$ is a normal crossing divisor in
$X_l$.

The action of each element $g\in G(k)$ lifts to the entire resolution sequence
\eqref{eqn5mar1}; this follows from~\cite[Theorem~13.2(2)(ii)]{bm}.
This means, inductively, that each blowup center $C_i$,
$i=0,\,1,\dots,\,l-2$, is invariant under this
action of $g$. Since we are assuming that $G(k)$ is Zariski dense in $G$
(see \S\ref{prel1.5}), each of these $C_i$ is
$G$-invariant; this implies that the action of $G$ lifts to the entire
resolution tower \eqref{eqn5mar1}, $C_{l-1}$ --- which is the strict transform
of $Y$ --- is $G$-invariant, the action of $G$ lifts to the blowup
\eqref{eqn13mar1}, and each $Y_i$, $i\le l$, is $G$-invariant.

In particular, $X_l$ is smooth, $Y_l$ is a $G$-invariant normal crossing divisor
in $X_l$ and the action of $G$ on $X_l-Y_l$ is free, since $Y_l$ contains the
preimage of $Y$. This implies that conditions (i) and (ii) 
of Definition~\ref{def3.1} are satisfied for $X_l$ and
the divisor $Y_l\subset X_l$.

We shall choose the centers $C_i$ for $i\ge l$ in such a way
that $C_i$ and $Y_{i}$ simultaneously have only normal crossings.
Inductively, this implies that for all $i\ge l$, $Y_i$ is a normal
crossing divisor and the action of $G$ on $X_i-Y_i$ is free.
With this choice of centers, 
conditions (i) and (ii) of Definition~\ref{def3.1} are satisfied for $X_i$ and
the divisor $Y_i$ for all $i\ge l$.

We would like the divisor $Y_n$ to satisfy
condition (iii) of Definition~\ref{def3.1}. In order to achieve this goal, 
we blow up, successively, all intersections of the components 
of the divisor $Y_l$, starting with those of the smallest dimension, as follows.

Let $m=\dim X$; then we define the center $C_{l}\subset X_l$ to be the
union of all $m$-tuple intersections of components of $Y_l$; it is a
finite set of points.
Inductively we define the center $C_{l+i}$ for
$i=1,\dots,m-2$ as the strict transform in $X_{l+i}$ of the union --- denote
it by $Y_l^{(i)}$ --- of $(m-i)$-tuple intersections of
components of $Y_l$.  Note that here $Y_l^{(i)}$ is a
union of smooth normal crossing $i$-dimensional subvarieties 
in $X_l$, and $C_{l+i}$ --- its strict
transform in $X_{l+i}$ --- is a union of disjoint smooth subvarieties;
similarly, the strict transform of $Y_l$ in $X_{l+m-1}$ is
the union of disjoint smooth subvarieties.
Each center $C_i$ we have described this far, is $G$-invariant, and
$C_i$ and $D_{i}$ simultaneously have only normal crossings.

Let $Z$ be an irreducible component of $Y_{l+m-1}$; it is
either (a) the strict transform of an irreducible component $Z'$ 
of $Y_l$ or (b) an irreducible
component of the exceptional divisor of the composition
\[ \pi_{l+m-1}\dots\pi_{l+1}\colon X_{l+m-1}\to X_l \; . \]

In case (a), for any $g\in G$ the subvariety $g(Z)$ is also the strict
transform of the irreducible component $g(Z')$ of $Y_l$; both $Z$ and
$g(Z)$ are components of the strict transform of $Y_l$ in $X_{l+m-1}$.
As the latter is the union of disjoint components,
either $g(Z)$ coincides with $Z$ or is disjoint from it. This means that
$Z$ satisfies condition (iii) of Definition~\ref{def3.1}.

In case (b) $Z$ satisfies condition (iii) of Definition~\ref{def3.1} by
Lemma~\ref{lem4.1-1}.

Therefore, the divisor $Y_{l+m-1}=D_{l+m-1}\cup\pi^{-1}(Y)$
satisfies condition (iii) of Definition~\ref{def3.1}, and consequently,
$X_n=X_{l+m-1}$ is in standard form with respect to it.
\end{pf*}

\begin{remark}
At the beginning of the proof of Theorem~\ref{thm.b1-2}, we could have 
taken an alternative approach by considering 
the canonical resolution of {\em the sheaf of ideals $\cI_Y$
of\/} $Y$ in $X$, as in \cite[Theorem~1.10]{bm}, instead of first considering
the canonical embedded resolution of singularities of $Y$, 
as in \cite[Theorem~1.6]{bm}, and then blowing up the strict 
transform $C_{l-1}$ of $Y$. Note that the action of $g\in G$ lifts 
to the canonical resolution of $\cI_Y$; this may be deduced from
\cite[Remark~1.5]{bm}. 

Alternatively, we could have used the constructive resolution of 
the idealistic space determined by the couple $(\cI_Y,1)$, 
as in \cite[Definition 2.4.1 and Theorem~7.3]{Vil2}.
The action of $g\in G$ lifts to this resolution by an argument
similar to that of~\cite[Corollary~7.6.3]{Vil2}.
\end{remark}

\begin{cor} \label{cor.b1-1}
Let $X$ be a $G$-variety and $Y\subset X$ a closed nowhere dense $G$-invariant
subvariety such that the action of $G$ on $X-Y$ is free. Then
there is a sequence of blowups
\begin{equation}  \label{eqn.cor.b1-1}
\pi\colon X_n \stackrel{\pi_n}{\lra} X_{n-1}
\dots \stackrel{\pi_2}{\lra} X_1 \stackrel{\pi_1}{\lra} X_0 = X
\end{equation}
where the centers $C_i\subset X_{i}$ are smooth and $G$-invariant, and $X_n$
is in standard form with respect to a divisor $\tilde Y\subset X_n$
which contains $\pi^{-1}(Y)$.
\end{cor}

\begin{pf}
Note that since $Y$ is nowhere dense in $X$, it is nowhere dense in each
irreducible component of $X$.

Consider the canonical resolution of singularities of $X$,
\begin{equation} \label{eqn14mar1}
X_l@>\pi_l>>\dots @>\pi_1>>X_0=X\ ,
\end{equation}
as in~\cite[Theorem~7.6.1]{Vil2} or~\cite[Theorem~13.2]{bm}.
The variety $X_l$ is smooth; similarly to the proof of Theorem~\ref{thm.b1-2},
we find that the centers $C_i\subset X_{i}$ are smooth 
and $G$-invariant, and the action of $G$ lifts to the entire resolution 
sequence \eqref{eqn14mar1}.

Let  $Y_l$ be the preimage of $Y$ in $X_l$. Then  $Y_l$ is nowhere dense in
each of the irreducible components of $X_l$, since $Y$ is nowhere dense in each
irreducible component of $X$. Consequently, $Y_l$ is nowhere dense in $X_l$.
Now apply Theorem~\ref{thm.b1-2} to $X_l$ and $Y_l$
to obtain a sequence
$X_n@>\pi_n>>\dots @>\pi_{l+1}>>X_l$
with smooth $G$-invariant centers, such that $X_n$ is in standard form
with respect to a divisor $\tilde Y\subset X_n$
which contains $\pi^{-1}(Y)$.
\end{pf}

%%%%%%%%%%%%%%%%%%%%%%%%%

\section{$G$-varieties in standard form}
\label{sect4}

With the exception of Remark~\ref{rem3.3c}, we shall assume 
throughout this section that the base field $k$ is algebraically closed.

\begin{thm} \label{thm3.2}
Let $X$ be a generically free $G$-variety in standard form,
and let $Y$ be as in Definition~\ref{def3.1}.
Suppose $x \in X$ lies on exactly $m$ irreducible components of $Y$. 
Then $\Stab(x)$ is isomorphic to a semidirect product $U \sdp A$, 
where $U$ is a unipotent group and $A$ is a diagonalizable
group of rank $\leq m$. 
\end{thm}

Our proof of Theorem~\ref{thm3.2} relies on the following lemma.

\begin{lem} \label{linearization}
Let $H$ be a diagonalizable group, $X$ an $H$-variety, and let
$X^H$ be the fixed point set of $H$ in $X$. If $X$ is smooth at a point $x$
and $x \in X^H$ then $X^H$ is also smooth at $x$; moreover, $T_x(X^H) =
T_x(X)^H$.
\end{lem}

\begin{pf} Note that if $X$ is affine then the lemma is a consequence
of the Luna Slice Theorem; see~\cite[Corollary to Theorem~6.4]{pv}.
Moreover, since every
quasiaffine $H$-variety can be equivariantly embedded into an affine  
$H$-variety (see~\cite[Theorem 1.6]{pv}), the lemma also holds if $X$ 
is quasiaffine. Thus it is sufficient to show that $x$ has 
an open quasiaffine $H$-invariant neighborhood $U \subset X$.

After replacing $X$ by its smooth locus (which is open and $H$-invariant),
we may assume $X$ is smooth.
Let $H^0$ be the identity component of $H$; since $H$ is diagonalizable,
$H^0$ is a torus (possibly $H^0=\{1\}$).
By a result of Sumihiro (see~\cite[Corollary~2]{sumihiro})
there exists an affine
$H^0$-invariant neighborhood $X_0$ of $x$ in $X$. We now define $U$ as
\[ U = \bigcap_{\overline{h} \in H/H^0} \overline{h}(X_0) \; .\]
Since $H/H^0$ is a finite group, $U$ is an open $H$-invariant quasiaffine
neighborhood of $x$, as claimed.
\end{pf}

\begin{pf*}{Proof of Theorem~\ref{thm3.2}}
Consider the Levi decomposition $\Stab(x) = U \sdp A$, where
$A$ is reductive and $U$ is unipotent; see, e.g., \cite[Section 6.4]{ov}.
We want to show that $A$ is, in fact, a diagonalizable group of rank $\leq m$. 

Denote the irreducible components of $Y$ passing through $x$ by
$Z_1, \dots, Z_m$; they intersect transversely at $x$.
Recall that by our assumption each $Z_i$ is $\Stab(x)$-invariant; hence,
their intersection $W=Z_1\cap\dots\cap Z_m$ is also $\Stab(x)$-invariant.

As $A$ is reductive, there is an $A$-invariant subspace $V$ in $T_x(X)$
complementary to $T_x(W)$.
We have an $A$-invariant decomposition
\[
V=V_1 \oplus V_2 \dots \oplus V_m
\]
where
\begin{equation}  \label{eqn.Vi}
V_i=V\cap T_x(Z_1)\cap\dots\cap\widehat{T_x(Z_i)}\cap\dots\cap T_x(Z_m)\ ;
\end{equation}
each $V_i$ is one-dimensional.
The group $A$ acts on each $V_i$ by a character, say,
$\chi_i\colon A\to \bbG_\ml$ (possibly trivial). 
We claim that the homomorphism
\[
\chi = (\chi_1, \dots, \chi_m)\colon A \lra (\bbG_\ml)^m
\]
is injective. Note that the theorem is an immediate consequence of this claim.

To prove the claim, note that $\Ker(\chi)$ is a reductive
subgroup of $A$. Thus in order to prove that $\Ker(\chi) = \{ 1\}$,
it is sufficient to show that every diagonalizable subgroup $\Ker(\chi)$ 
is trivial.  (Indeed, this immediately implies that the identity component 
$\Ker(\chi)^0$ is unipotent and, hence, trivial; 
see~\cite[Exercise 1, p. 137]{humphreys}. Thus $\Ker(\chi)$ is finite 
and every abelian subgroup of $\Ker(\chi)$ is trivial; this is only 
possible if $\Ker(\chi) = \{ 1 \}$.)  

Let $H \subset \Ker(\chi)$ be a diagonalizable group; we want 
to show that $H = \{ 1 \}$. Assume the contrary.
Denote the fixed point set of $H$ by $X^H$. Since the action of $G$ on
$X-Y$ is free, $X^H \subset Y$.
By Lemma~\ref{linearization}, 
$X^H$ is smooth. Consequently, only one irreducible component of $X^H$ 
passes through $x$; denote this component by $X^H_0$. Then $X^H_0$ is 
contained in one of the components $Z_1, \dots, Z_m$, say in $Z_i$, and
by Lemma~\ref{linearization}, 
\begin{equation} \label{e3.1}
T_x(X)^H = T_x(X^H) = T_x(X^H_0) \subset T_x(Z_i)  \; .
\end{equation}
Now note that by our assumption, $\chi_i|_H$ is trivial and thus
$V_i \subset T_x(X)^H$ but, on the other hand, by \eqref{eqn.Vi}
$V_i \not \subset T_x(Z_i)$, 
contradicting \eqref{e3.1}. This completes the proof of the claim.
\end{pf*}

\begin{remark} \label{rem3.3a}
Note the following interesting special cases of Theorem~\ref{thm3.2}:

\smallskip
(a) If $\Stab(x)$ is connected then $\Stab(x)$ is solvable and

\smallskip
(b) if $\Stab(x)$ is finite then $\Stab(x)$ is commutative.
\end{remark}

\begin{remark} \label{rem3.3b}
Our proof shows that 
\begin{equation} \label{eqn.rem.31}
T_x(X)/T_x(W)=\bigoplus_{i=1}^m
\frac{T_x(Z_1)\cap\dots\cap\widehat{T_x(Z_i)}\cap\dots\cap T_x(Z_m)}{T_x(W)}\ .
\end{equation}
is a direct sum decomposition of the normal space
$T_x(X)/T_x(W)\isomo V$ as a direct sum of 1-dimensional character spaces
for the natural action of $A$. Moreover, the above (diagonal)
representation of $A$ on $V$ is faithful.
\end{remark}

\begin{remark} \label{rem3.3c}
Suppose the base field $k$ is not necessarily algebraically closed
(but is of characteristic 0),
$X$ is a generically free $G$-variety in standard form,
and $x \in X$
has a finite stabilizer of exponent $e$. Then the residue field $k'$ of $x$
contains a primitive $e$-th
root of unity. Indeed, $\Stab(x)$ has a faithful diagonal 
representation~\eqref{eqn.rem.31} defined over $k'$; this is only possible
if $k'$ contains a primitive $e$-th root of unity.
\end{remark}

\begin{cor} \label{cor3.3}
Let $X$ be a generically free $G$-variety 
in standard form.  Suppose that $H= \Stab(x)$ is a finite group.
Then $\dim(X) \geq \dim(G) + \rank(H)$.
\end{cor}

Here $\rank(H)$ denotes the rank of the finite abelian group
$H= \Stab(x)$; see Remark~\ref{rem3.3a}(b).

\begin{pf}
Let $Y$ be as in Definition~\ref{def3.1}. Suppose
exactly $m$ irreducible components $Z_1, \dots, Z_m$ meet at $x$; then by
Theorem~\ref{thm3.2} we have $m \geq \rank(H)$.
Since $Z_1, \ldots, Z_m$ intersect transversely at $x$, their intersection
$W = Z_1 \cap \dots \cap Z_m$ is smooth at $x$ and
\[
 \dim(X) = \dim_x(W) + m = \dim(W_0) + m \; ,
\]
where $W_0$ is the (unique) component of $W$ passing through $X$. 

Since $m \geq \rank(H)$,
it only remains to show that $\dim(W_0) \geq \dim(G)$. Indeed, let
\[
G' = \{ g \in G \mid g(Z_i) = Z_i  \; \forall i = 1, \dots, m \} \; .
\] 
Then $G'x \subset W_0$. Since $G'$ is a subgroup of finite index in $G$
and $\Stab(x)$ is assumed to be finite, 
we have $\dim(W_0) \geq \dim(G'x) = \dim(G') = \dim(G)$,
as claimed.
\end{pf}

%%%%%%%%%%%%%%%%%%%%%%%%%%%

\section{The behavior of fixed points under rational morphisms}
\label{sect5}

Suppose $H$ is an algebraic group and $f\colon X \brokrarr Y$ is a rational 
map of $H$-varieties. In this section we shall be interested in two types 
of results (under certain additional assumptions on $H$, 
$X$, $Y$ and $f$): {\em ``going down"} results, which assert that 
if $H$ fixes a point of $X$ then it fixes a point of $Y$ and 
{\em ``going up"} results which assert the converse. 

Note that the ``going down" assertion is always true
if $f$ is a regular map; indeed, if $x \in X$ is fixed by $H$ then so is
$f(x) \in Y$. The situation is somewhat more complicated for rational maps; 
in particular, we need to make a strong assumption on the group $H$; see 
Remark~\ref{rem.Going-down}.

Throughout this section we shall assume that the base field $k$ is 
algebraically closed.

The proofs we originally had in this section relied on canonical resolution 
of singularities; cf.\ Remark~\ref{rem6.2}. 
Koll\'ar and Szab\'o recently found simple characteristic-free 
proofs of Propositions~\ref{prop6.2} and~\ref{prop6.3}. 
These proofs are presented in the Appendix at the end of this paper; we shall 
therefore omit most of our original arguments. We also note that 
our earlier versions of Lemma~\ref{lem4.1} 
and Proposition~\ref{prop6.2} assumed that $H$ is diagonalizable; our
earlier version of Proposition~\ref{prop6.2} (respectively, 
Proposition~\ref{prop6.3}) assumed that $Y$ (respectively 
$X$) is projective, rather than complete. The current 
Propositions~\ref{prop6.2} and~\ref{prop6.3} are characteristic zero versions
of, respectively, Propositions~\ref{Going-down} and~\ref{Going-up}.

We begin with a simple lemma.

\begin{lem} \label{lem4.1}
Let $H = U \sdp A$, where $U$ is unipotent and $A$ is diagonalizable,
$X$ be an $H$-variety and
$\pi\colon X_1 \lra X$ be a blowup with a smooth $H$-invariant center 
$C \subset X$.
If $x$ is a smooth point of $X$ which is fixed by $H$ then there exists
an $x_1 \in X_1$ such that $\pi(x_1) = x$ and $x_1$ is fixed by $H$.
\end{lem}

\begin{pf} Recall that $\pi$ is 
an isomorphism over $X-C$; thus if $x \not \in C$ then we can take
$x_1 = \pi_1^{-1}(x)$. On the other hand, if $x \in C$ then $f^{-1}(x)
\simeq \bbP(V)$ (as $H$-varieties), where $V = N_x(C) = T_x(X)/T_x(C)$. 
The action of $H$ has an eigenvector in $V$ (see Lemma~\ref{class.defn});
thus $H$ fixes some $x_1 \in\bbP(V) = f^{-1}(x)$, as claimed.
\end{pf}

\begin{remark} \label{rem.thm1.1-optimal}
Lemma~\ref{lem4.1} shows that Theorem~\ref{thm1.1} is sharp in the sense that
the stabilizers of points of $X_n$ cannot be further reduced by additional 
blowups with smooth equivariant centers. 
\end{remark}

\subsection*{Going down}

\begin{prop} \label{prop6.2} 
Let $H = U \sdp A$, where $U$ is unipotent and $A$ is diagonalizable.
Suppose $f\colon X \brokrarr Y$ is a dominant rational map of 
$H$-varieties, where $Y$ is complete. 
If $H$ fixes a smooth point $x$ in $X$ then $H$ fixes a point $y \in Y$.
\end{prop}

\begin{pf} See Proposition~\ref{Going-down}.
\end{pf}

\begin{remark} \label{rem6.2}
We will now briefly outline our original proof of Proposition~\ref{prop6.2}.
It is more complicated than the proof of Proposition~\ref{Going-down}
and only works in characteristic zero; however, we feel this argument 
may be of independent interest.

First we showed that there exists a sequence of blowups
\[ X_n \stackrel{\pi_n}{\lra} \dots \stackrel{\pi_2}{\lra} X_1 
 \stackrel{\pi_1}{\lra} X_0 = X \]
with smooth $H$-invariant centers such that $f$ lifts to a regular map 
\[ f'\colon X_n \lra Y \]
of $H$-varieties. This is, in fact, true for any algebraic group $H$ and any
$H$-invariant rational map $f \colon X \brokrarr Y$; the proof relies on
canonical resolution of singularities.

Applying Lemma~\ref{lem4.1} inductively to the above tower of blowups, 
we see that for every $i = 0, 1, \dots, n$ there exists a (necessarily smooth)
$H$-fixed point $x_i\in X_i$ lying above $x=x_0$. 
Now $y = f'(x_n)$ is an $H$-fixed point of $Y$.
\qed
\end{remark}

\subsection*{Going up}

Let $H$ be a diagonalizable group,
$f\colon X \brokrarr Y$ be a rational map of $H$-varieties. We now want to
prove that if $H$ fixes a smooth point $y \in Y$ then $H$ fixes a point
of $X$. We clearly need to assume that $f$ is dominant
and the fibers of $f$ are complete; the following example shows that
these assumptions are not sufficient, even if $X$ is irreducible.

\begin{example}
Let $H = \bbZ/n_1\bbZ \times \dots \times \bbZ/n_r\bbZ $ 
be a finite abelian group, $Y$ be an $H$-variety, 
$P$ be a projective $H$-variety where $H$ acts freely (i.e., all stabilizers
are trivial), and $X = Y \times P$. Then $H$ acts freely on $X$,
hence, the ``going up" assertion will fail for the map $f \colon X \lra Y$,
where $f$ = projection to the first component.
(Note that the fibers of this map are projective, so lack of completeness is
not the problem here.) 
To construct $P$, let $E$ be an elliptic curve and let 
$p_i$ be a point of order $n_i$ on $E$. Now set $P = E^n$ 
and define the $H$-action on $P$ by
\[ (i_1, \dots, i_r) \cdot (x_1, \ldots, x_n) = (x_1 + i_1p_1,
\dots, x_r + i_r p_r) \; , \]
where $+$ refers to addition on $E$.
\qed
\end{example}

Nevertheless, it turns out that one can still prove 
a useful ``going up" property.  

\begin{prop} \label{prop6.3} 
Let $H$ be an abelian $p$-group 
and $f\colon X \brokrarr Y$ be a dominant rational $d:1$-map
of generically free $H$-varieties. Assume $X$ is complete,
$d$ is prime to $p$ and $y \in Y$ is a smooth point 
fixed by $H$. Then $H$ fixes a point $x \in X$.
\end{prop}

\begin{pf}
Note that since $y$ is a smooth point of $Y$ fixed by $H$, the
irreducible component $Y_0$ of $Y$ containing $y$, is preserved by $H$.
Replacing $Y$ by $Y_0$ and $X$ by the union of its
irreducible components which are mapped dominantly onto  $Y_0$, we may
assume that $Y$ is irreducible and each component $X_i$ of $X$ is mapped
dominantly onto $Y$.

Similarly to the argument of the proof of Proposition~\ref{Going-up},
we note that $H$ acts on the set $\{X_i\}$;
let $\cX_j$ be the $H$-orbits in this set.
Pick an element $X_j^*$ in $\cX_j$; then
\[
d=\deg(X/Y)=\sum_j|\cX_j|\cdot\deg(X_j^*/Y)\ .
\]
As $d$ is not divisible by $p$, there is an orbit $\cX_0$ consisting
of a single element $X_0^*$ such that $\deg(X_0^*)$ is not divisible
by $p$.
Replacing $X$ by $X_0^*$, we may assume that $X$ is irreducible; now
apply Proposition~\ref{Going-up}.
\end{pf} 

\section{Essential dimensions and cohomological invariants}

\subsection*{Essential dimension}

We now recall the definition of essential dimension from \cite{r2};
in the case of finite groups, see also~\cite{br1} and \cite{br2}.

\begin{defn} \label{def.ed}
(1) The essential dimension of a primitive generically free $G$-variety 
$X$ is the minimal value of $\dim(Y/G) = \dim(Y) - \dim(G)$, where
$Y/G$ denotes the rational quotient of $Y$ by $G$ 
and the minimum is taken over all $G$-compressions $X \brokrarr Y$;
see~\S\ref{prel2} and~\S\ref{prel4}. We denote this number by $\ed(X)$.  

\smallskip
(2) If $V$ is a generically free irreducible linear 
representation of $G$, we refer to $\ed(V)$ as the essential 
dimension of $G$ and denote it by $\ed(G)$. By~\cite[Theorem~3.4]{r2} 
this number is independent of the choice of $V$. 
Equivalently, $\ed(G)$ can be defined as the maximal value of $\ed(X)$, 
as $X$ ranges over all primitive generically free $G$-varieties; 
see~\cite[Section~3.2]{r2}.
\end{defn}

\begin{remark} \label{rem6.05} The definition of essential dimension 
of an algebraic group
in~\cite{r2} assumes that the base field $k$ is algebraically closed
and of characteristic 0;
the definition of essential dimension of a finite group in~\cite{br1} 
and \cite{br2} is valid over an arbitrary field of characteristic 0.
In this paper we will be interested, almost exclusively,
in proving lower bounds on essential dimensions of various groups 
and $G$-varieties.  Since $\ed(X) \geq \ed(X \otimes_k \kbar)$ 
for any $G$-variety $X$, with $G$ finite, as well as for any imaginable 
notion of $\ed(X)$ with $G$ infinite, a lower bound on $\ed(G)$ or
$\ed(X)$ over $\kbar$ will automatically be valid over $k$.
For this reason, all lower bounds we prove under the assumption that $k$ is
algebraically closed, also hold without this assumption.
\end{remark}

\subsection*{Essential dimension at $p$}

We will also study the following related numerical invariants which
were brought to our attention by J.-P. Serre. 

\begin{defn} \label{def.edp} 
(1) Let $p$ be a prime integer and let $X$ be a primitive generically free
$G$-variety.  We define the {\em essential dimension of X at} $p$ 
as the minimal value of $\ed(X')$, where the minimum is taken 
over all dominant rational $d:1$ maps $X' \brokrarr X$ of primitive
$G$-varieties (see~\S\S\ref{prel1}, \ref{prel2} and \ref{prel4}),
with $d$ prime to $p$.  We shall denote this number by $\ed(X; p)$.

(2) The {\em essential dimension of $G$ at $p$} is defined as the
maximal value of $\ed(X; p)$, as $X$ ranges over all primitive
generically free $G$-varieties.   We shall denote this number by $\ed(G; p)$.
\end{defn}

\begin{remark} \label{rem.edp} $\ed(X; p)$ is closely related to
the ``relative essential dimension" $\ed^{m, H}(X; p)$
defined (for finite groups only) in \cite[Section 5]{br2}.
More precisely,
$\ed(X; p) $ is the maximal value of $\ed^{m, H}(X; p)$, as $H$ ranges over 
all finite groups and $m$ ranges over all positive integers prime to $p$.
We shall not work with $\ed^{H, m}(X)$ in this paper.
\end{remark}

\begin{remark} \label{rem.except-prime}
Clearly, $\ed(X) \geq \ed(X; p)$ for every primitive generically free
$G$-variety $X$ and every prime $p$.
In particular, $\ed(G) \geq \ed(G; p)$. Note also that if $G$ is a simple
group then $\ed(X; p)=0$ unless $p$ is one of the so-called exceptional 
primes. For details, including a list of exceptional primes, 
see~\cite[Section 2]{serrepp}. 
\end{remark}

The following lemma will not be needed in the sequel; we include it here to
illustrate the similarity between the definitions of 
$\ed(G)$ and $\ed(G; p)$. 

\begin{lem} \label{lem.ed-p}
Suppose $G$ is an algebraic group and $p$ is a prime integer. 

\smallskip
(a) Let $X$ be a primitive generically free $G$-variety and
$f\colon X \brokrarr Y$ be a $G$-compression. Then $\ed(X; p) \leq \ed(Y; p)$. 

\smallskip
(b) Let $V$ be a generically free linear representation of $G$.
Then $\ed(V; p) = \ed(G; p)$. In other words, $\ed(V; p) \geq
\ed(X; p)$ for any primitive generically free $G$-variety $X$; 
in particular, $\ed(V; p)$ is independent of the choice of $V$.
\end{lem}

\begin{pf} 
(a) Suppose $Y' \brokrarr Y$ is a $d:1$ dominant rational 
map of primitive $G$-varieties. It is enough to show that there exists
a commutative diagram of rational maps 
\begin{equation} \label{e6.6} 
\begin{array}{rcl}
           X' &  \stackrel{f'}{\brokrarr} &  Y' \\
           | \; &         & \; |   \\
    e:1 \;   |  \; &         & \; | \; d:1  \\
  \downarrow \, &         & \, \downarrow   \\
           X  & \stackrel{f}{\brokrarr} &  Y 
\end{array}
\end{equation}
of primitive $G$-varieties, where $X' \brokrarr X$ 
is an $e:1$ dominant rational
map of primitive $G$-varieties and $e$ is not divisible by $p$. Indeed, 
the existence of $f'$ immediately implies
$\ed(X') \leq \ed(Y')$ (cf. \cite[Lemma 3.3(b)]{r2});
taking the minimum over all $Y'$, we obtain the desired inequality.

To construct the diagram \eqref{e6.6}, note that since $X$, $Y$ and $Y'$
are primitive, $k(X)^G$, $k(Y)^G$ and $k(Y')^G$ are, by definition, fields; 
see~\S\ref{prel4}.
Moreover, $[k(Y')^G:k(Y)^G] = d$. We claim that there exists a diagram
of field extensions

\setlength{\unitlength}{.7pt}
\begin{picture}(320,175)
\put(280,10){$k(Y)^G$}
\put(135,57){$k(X)^G$}
\put(280,95){$k(Y')^G$}
\put(145,142){$L$}

\put(138,102){$e$}
\put(301,55){$d$}

\put(271,22){\line(-3,1){95}}
\put(271,107){\line(-3,1){95}}
\put(295, 33){\line(0,1){50}}
\put(150,80){\line(0, 1){50}}
\end{picture} 

\noindent
where $L$ contains both $k(X)^G$ and $k(Y')^G$ and $p \notdiv e= [L:k(X)^G]$.
Indeed, write 
\[ k(X)^G \otimes_{k(Y)^G} k(Y')^G =
L_1 \oplus \dots \oplus L_m \; , \]
where each $L_i$ is a field;
see~\cite[Section 5.6]{jacobson.rt}. Since  
\[ \sum_{i=1}^m [L_i: k(X)^G] = \dim_{k(X)^G} \, 
\Bigl(k(X)^G \otimes_{k(Y)^G} k(Y')^G\Bigr) =[k(Y')^G: k(Y)^G] = d \]
is not divisible by $p$, we conclude that
$p \notdiv \, [L_i: k(X)^G]$ for some $i$. Now set $L=L_i$ and
$e = [L_i:k(X)^G]$.

The above diagram gives rise to the following diagram of rational maps:
\[ \begin{array}{rcl}
           X_0  & \stackrel{\overline{f'}}{\brokrarr} &  Y'/G \\ 
           | \; &         & \; |   \\
    e:1 \;   |  \; &         & \; | \; d:1  \\
  \downarrow \, &         & \, \downarrow   \\
           X/G  & \stackrel{\overline{f}}{\brokrarr} &  Y/G 
\end{array} \]
where $X_0$ is an irreducible algebraic variety whose function field in $L$. 
Taking the fiber product of this diagram with $Y$ over $Y/G$, and remembering
that $Y/G \times_{Y/G} Y \simeq Y$, $Y'/G \times_{Y/G} Y \simeq Y'$,
and $X/G \times_{Y/G} Y \simeq X$ as $G$-varieties 
(see~\cite[Lemma 2.14]{r2}), we obtain the desired diagram \eqref{e6.6} with 
$X' = X_0 \times_{Y/G} Y$.
Note that $X'/G \simeq X_0$ and thus $X'$ is a primitive $G$-variety;
see~\cite[Lemma 2.14]{r2}.

(b) Recall that by \cite[Corollary 2.17]{r2}, for every primitive
generically free
$G$-variety $X$, there exists a $G$-compression $X \times \bbA^d \brokrarr V$,
where $d = \dim(V)$ and $G$ acts trivially on $\bbA^d$. (This fact is 
a consequence of the ``no-name lemma".) Thus by part (a) 
 \begin{equation} \label{e5.5}
\ed(X \times \bbA^d; p) \leq \ed(V; p)  \; .
\end{equation}
On the other hand, the argument of \cite[Lemma 5.3]{br2} shows that 
\[ \ed(X \times \bbA^1; p) = \ed(X;p) \]
for any primitive generically free $G$-variety $X$; see~Remark~\ref{rem.edp}.
This, along with \eqref{e5.5}, proves part (b). 
\end{pf}

\subsection*{Cohomological invariants}

A simple but important relationship between the essential dimension 
of an algebraic group $G$ and its cohomological invariants was 
observed by J.-P. Serre (see Lemma~\ref{lem.ci} below). This observation
makes it possible to deduce lower bounds on $\ed(G; p)$ from the existence
of non-trivial cohomological invariants.

In the next section we will develop a method for proving
lower bounds on $\ed(G; p)$, which does not presuppose 
the existence of a non-trivial cohomological
invariant.  However, for the purpose of motivating our results and
placing them in the proper context, we briefly explain 
the relationship between cohomological invariants 
and essential dimension. We will follow up on this theme in
Remark~\ref{rem.ci2}.

Suppose $F$ is field, $\overline{F}$ is the algebraic closure of $F$, 
$\Gamma = \galois(\overline{F}, F)$ and $M$ 
is a torsion $\Gamma$-module. In the sequel, we shall denote the
Galois cohomology group by $H^i(F, M)$; see \cite{serregc}. 

We shall view
$H^i(\, \cdot \, , M)$ as a functor from the category of fields to 
the category of groups. We shall also consider the functor
$H^1( \, \cdot \, , G)$ from the category of finitely generated field
extensions of $k$ to the category of sets. Recall that elements
of the non-abelian cohomology set $H^1(F, G)$ are 
in 1---1 correspondence with primitive
generically free $G$-varieties $X$ such that $k(X)^G = F$;
see~\cite[I.5.2]{serregc},~\cite[Theorem 1.3.3]{popov} 
or~\cite[Lemma 12.3]{r2}.

\begin{defn} \label{def.ci}
A cohomological invariant $\alpha$ of $G$-varieties
is a morphism of functors $H^1(\, \cdot \,, G) \lra H^d(\, \cdot \, , M)$. 
In other words, $\alpha$ assigns a cohomology class 
$\alpha(X) \in H^d(k(X)^G, M)$ to every primitive generically 
free $G$-variety $X$, so that for every compression $X \brokrarr Y$, 
$\alpha(X)$ is the image of $\alpha(Y)$ under the natural restriction 
homomorphism $H^d(k(Y)^G, M) \lra H^d(k(X)^G, M)$. 
\end{defn}

\begin{remark} \label{rem.ci1} The above notion of cohomological invariant
(and the equivalent notion used in \cite[Section 12]{r2})
are somewhat more narrow than the usual definition (see \cite[6.1]{serrepp}
or \cite[31B]{boi}), due to the fact that we work over an algebraically 
closed field $k$. This means
that a cohomological invariant in the sense of \cite[Section 6.2]{serrepp}
or \cite[Section 31B]{boi} is also a cohomological
invariant in our sense but the converse may not be true.
\end{remark}

The following observation, due to J.-P. Serre, relates
the essential dimension $G$ to cohomological invariants.

\begin{lem} \label{lem.ci} 
Let $G$ be an algebraic group. Suppose there exists 
a non-trivial cohomological invariant $\alpha\colon H^1(\, \cdot \,, G) \lra 
H^i(\, \cdot \, , M)$, where $M$ is a $p$-torsion module. Then
$\ed(G; p) \geq i$. 
\end{lem}

\begin{pf} It is sufficient to show that if $\ed(G; p) < i$ then
that $\alpha(X) =0$ for every generically free primitive $G$-variety $X$. 

Indeed, for every generically free primitive $G$-variety $X$ 
there exists a $d:1$-cover $X' \brokrarr X$ of $G$-varieties and 
a $G$-compression $X' \lra Y$ such that $\trdeg_k \, k(Y)^G = \dim(Y/G) < i$. 
Thus $H^i(k(Y)^G, M) = (0)$ (see \cite[II.4.2]{serregc})
and consequently $\alpha(Y) = 0$. Since $\alpha(X')$ is a homomorphic
image of $\alpha(Y)$, we conclude $\alpha(X') = 0$. Finally, since
$[k(X')^G : k(X)^G] = d$ is prime to $p$, the restriction map
$H^i(k(X)^G, M) \lra H^i(k(X')^G, M)$ is injective; 
see~\cite[I.2.4]{serregc}.  Thus $\alpha(X) = 0$, as claimed.
\end{pf}

\section{Stabilizers as obstructions to compressions}
\label{sect7}

In this section we assume that $k$ is algebraically closed; see 
Remark~\ref{rem6.05}. 

\subsection*{A lower bound}

We begin by recalling the following result of Sumihiro.

\begin{prop} \label{prop.s-k}
Every $G$-variety is birationally isomorphic to a complete $G$-variety.
\end{prop}

\begin{pf} Let $X$ be a $G$-variety. 
After removing the singular locus from $X$, we may assume that
$X$ is smooth.  Then $X$ is a disjoint union of smooth irreducible 
varieties.  The group $G$ acts on the set of
irreducible components of $X$; the orbits of this
action give a decomposition of $X$ as 
a disjoint union of primitive $G$-varieties; cf.~\cite[Lemma~2.2(a)]{r2}.
Thus we may assume $X$ to be smooth and primitive. 

Let $X_0$ be an irreducible component of $X$, let $G_0$
be the subgroup of $G$ that preserves $X_0$, and let $X_0'$ be Sumihiro's
equivariant completion of $X_0$ as an irreducible $G_0$-variety; see
\cite[Theorem 3]{sumihiro}.
Then $X'=X_0'\times_{G_0}G$ is a $G$-equivariant completion of $X$;
it is a disjoint finite union of copies of $X_0'$.
In particular, $X'$ and $X$ are birationally isomorphic as $G$-varieties.
\end{pf}

We are now ready to prove our first lower bound on the essential dimension 
of a $G$-variety.

\begin{thm} \label{thm6.1} 
Let $H$ be a finite abelian subgroup of $G$ such that
\begin{itemize}
\item[(a)] the centralizer of $H$ is finite, and
\item[(b)] $H$ does not normalize any non-trivial unipotent subgroup of $G$.
\end{itemize}
Suppose $X$ is a primitive generically
free $G$-variety, $x$ is a smooth point of $X$ fixed by $H$, and
$X \brokrarr Y$ is a $G$-compression. Then
\begin{enumerate}
\item $\dim(X) \geq \rank(H) + \dim(G)$.
\item Moreover, 
$\dim(Y) \geq \rank(H) + \dim(G)$. In other words, 
$\ed(X) \geq \rank(H)$.
\item If $H$ is a $p$-group then $\ed(X; p) \geq \rank(H)$.
\end{enumerate}
\end{thm}

Note that since $X$ is primitive, $\dim(X)$ is the dimension of every 
irreducible component of $X$.  Moreover, since $X$ is primitive, so is $Y$; 
hence, $\dim(Y)$ is the dimension of
every irreducible component of $Y$.

\smallskip
\begin{pf} 
(1) By Corollary~\ref{cor.b1-1} there exists a tower
\[
\pi\colon X_n \stackrel{\pi_n}{\lra} X_{n-1}
\dots \stackrel{\pi_2}{\lra} X_1 \stackrel{\pi_1}{\lra} X_0 = X
\]
of blowups with smooth $G$-invariant centers such that $X_n$ is 
in standard form. Thus, in view of Lemma~\ref{lem4.1}, we may replace  
$X$ by $X_n$, i.e., we may assume without loss of generality that $X$
is in standard form.

By Theorem~\ref{thm3.2} $\Stab(x) = U \sdp A$, 
where $U$ is unipotent and $A$ is diagonalizable. Recall that 
$H \subset \Stab(x)$.  Since $U$ is normal
in $\Stab(x)$, it is normalized by $H$. Hence, in view of assumption (b),
we conclude that $U = \{ 1 \}$ and thus $\Stab(x) = A$. In particular,
$A \subset C_G(H)$; thus $A$ is finite. Now by Corollary~\ref{cor3.3} 
$$\dim(X) \geq \rank(A) + \dim(G) \geq \rank(H) + \dim(G) \; , $$
as claimed.

\smallskip
(2) By Proposition~\ref{prop.s-k} we may assume $Y$ is complete.  
Moreover, in view 
of Corollary~\ref{cor.b1-1} we may also assume that $Y$ is smooth.
By Proposition~\ref{prop6.2}, there exists a point $y \in Y$ such 
that $H \subset \Stab(y)$. Now apply part (1) to $Y$.

\smallskip
(3) Let $X' \brokrarr X$ be a $G$-equivariant $d:1$-cover of $X$.
We want to show $\ed(X') \geq \rank(H)$.
By Proposition~\ref{prop.s-k} we may assume $X'$ is complete; moreover, 
by Corollary~\ref{cor.b1-1} we may also assume $X'$ is smooth. 
By Proposition~\ref{prop6.3} there exists a point
$x' \in X'$ that is fixed by $H$. We now apply part (2) to 
$X'$ to conclude that $\ed(X') \geq \rank(H)$.
\end{pf}

\begin{cor} \label{cor5.5}
Let $G$ be an algebraic group and $H$ be an abelian subgroup 
of $G$ such that
(a) the centralizer of $H$ is finite and (b) $H$ does not normalize
any non-trivial unipotent subgroup of $G$.
Then $\ed(G) \geq \rank(H)$. Moreover, if $H$ is a $p$-group
then $\ed(G, p) \geq \rank(H)$. 
\end{cor}

\begin{pf} Apply Theorem~\ref{thm6.1}(2) and (3)
to $X = V$ = generically free 
linear representation of $G$ and $x = 0 \in V$.
\end{pf}

\begin{example} Let $G$ be a finite group and $H \simeq (\bbZ/p\bbZ)^m$
be a subgroup of $G$. Then $\ed(G; p) \geq m$. In particular, $\ed(\Sym_n; p)
\geq [n/p]$; cf. \cite[Section~6.1]{br1} and \cite[Section~7]{br2}.
\end{example}

\subsection*{A lemma of Serre}

The difficulty in applying Theorem~\ref{thm6.1} 
and Corollary~\ref{cor5.5} is that
condition (b) is often hard to verify. Fortunately,
under rather general assumptions, there is an easy way around 
this problem.

\begin{remark} Let $G$ be an algebraic group. Assume there exists
an abelian subgroup $H$ of $G$ satisfying 
conditions (a) and (b) of Theorem~\ref{thm6.1}. Then the identity 
component of $G$ is semisimple.
\end{remark}

\begin{pf} Assume $G$ is not reductive. Then the unipotent 
radical $R_u(G)$ is a non-trivial normal unipotent subgroup
of $G$, and thus condition (b) fails.  

Now assume 
$G$ is reductive. The radical $R(G)$ is the connected component 
of the center of $G$ (see~\cite[19.5]{humphreys}); hence, condition (a)
fails unless $R(G)$ is trivial. This means that the identity component 
of $G$ is semisimple, as claimed.
\end{pf} 

Thus if $G$ is connected, we may assume without loss of generality
that it is semisimple. The following lemma, communicated to us by
J.-P. Serre, shows that in this case conditions (a) and (b) of
Theorem~\ref{thm6.1} are equivalent.

\begin{lem} \label{lem6.11}
Let $G$ be a connected semisimple group and let $H$ be 
a (not necessarily connected) reductive subgroup of $G$.
Then the following conditions are equivalent.

\smallskip
(a) The centralizer $C_G(H)$ of $H$ in $G$ is infinite.

\smallskip
(b) $H$ normalizes a non-trivial unipotent subgroup of $G$.

\smallskip
(c)  $H$ is contained in a proper parabolic subgroup of $G$.
\end{lem}

\begin{pf} We will first show that (c) $\Longrightarrow$ (b),
then use this implication to prove that
(a) $\Longrightarrow$ (b) $\Longrightarrow$ (c) $\Longrightarrow$ (a).

\smallskip
(c) $\Longrightarrow$ (b): If $H$ is contained in a proper parabolic
subgroup $P$ then $H$ normalizes the unipotent radical $R_u(P) \neq \{ 1 \}$. 

\smallskip
(a) $\Longrightarrow$ (b): Assume $C_G(H)$ is infinite. If $C_G(H)$ 
contains a non-trivial
unipotent element $u$ then $H$ centralizes (and, hence, normalizes) 
the unipotent subgroup $\overline{\lf<u\r>} \neq \{ 1 \}$ and thus (b) holds.
If the centralizer $C_G(H)$ does not contain a non-trivial unipotent element,
then the identity component of $C_G(H)$ is a non-trivial torus $T$. 
In this case $H \subset C_G(T)$, and $C_G(T)$ is a Levi subgroup of some
non-trivial parabolic subgroup of $G$; see~\cite[30.2]{humphreys}.
Thus (c) holds, and, hence, so does (b).

\smallskip
(b) $\Longrightarrow$ (c): Suppose $H$ normalizes a non-trivial 
unipotent subgroup $U$ of $G$. Recall that the Borel---Tits 
construction associates, in a canonical way, a parabolic 
subgroup $P(U)$ to $U$ so that $U$ is contained in 
the unipotent radical of $P(U)$; see \cite[30.3]{humphreys}. 
In particular, $P(U)$ is proper.  Moreover, by our 
assumption $H \subset N_G(U)$, where $N_G(U)$ denotes 
the normalizer of $U$ in $G$. Since $N_G(U) \subset P(U)$
(see \cite[Corollary 30.3A]{humphreys}), $H$ is contained in the
proper parabolic subgroup $P(U)$. This proves (c).

\smallskip
(c) $\Longrightarrow$ (a): If $H$ is contained in a proper 
parabolic subgroup 
$P$ of $G$ then, by Levi's decomposition theorem, $H$ is contained 
in some Levi subgroup $L$ of $P$; see~\cite[Theorem 6.4.5]{ov}.
Then $C_G(L) \subset C_G(H)$. Since the center $Z(L)$ contains 
a non-trivial torus (see~\cite[30.2]{humphreys}), and $Z(L) \subset C_G(L)
\subset C_G(H)$, we conclude that $C_G(H)$ is infinite. 
\end{pf}

\subsection*{A better bound}

We can now prove the main results of this section.

\begin{thm} \label{thm5-5}
Let $G$ be an algebraic group, $H$ be an abelian subgroup of $G$,
and $X$ is a generically free $G$-variety. Suppose $H \subset \Stab(x)$ for
a smooth point $x$ of $X$.

\smallskip
(1) Assume $G$ is (connected and) semisimple and the centralizer $C_G(H)$ 
is finite. Then $\ed(X) \geq \rank(H)$. Moreover, if $H$ is a $p$-group then
$\ed(X; p) \geq \rank(H)$. 

\smallskip
(2) More generally, if the identity component $G^0$ of $G$ is semisimple
and the centralizer $C_{G^0}(H \cap G^0)$ is finite then $\ed(X) \geq 
\rank(H)$. Moreover, if $H$ is a $p$-group then
$\ed(X; p) \geq \rank (H)$. 
\end{thm}

\begin{pf} It is enough to verify that $G$ and $H$ satisfy conditions
(a) and (b) of Theorem~\ref{thm6.1}. In part (1) this follows immediately 
from Lemma~\ref{lem6.11}. 

(2) To check condition (a), note that $C_{G^0}(H \cap G^0)$ is of 
finite index in $C_{G}(H \cap G^0)$. This implies that $C_G(H \cap G^0)$ 
is finite and, hence, so is $C_G(H)$. To check condition (b),
note that since we are working over a field of characteristic 0, unipotent 
subgroups of $G$ are connected (see, e.g.,~\cite[3.2.2, Corollary 2]{ov}) 
and, hence, contained in $G^0$. 
By Lemma~\ref{lem6.11}, $H \cap G^0$ does not normalize any of 
them (except for $\{ 1 \}$). Hence, neither does $H$.
\end{pf}

\begin{thm} \label{thm5-6}
Let $G$ be an algebraic group and $H$ be an abelian subgroup of $G$.

\smallskip
(1) Suppose $G$ is (connected and) semisimple and the centralizer $C_G(H)$ 
is finite. Then $\ed(G) \geq \rank(H)$. Moreover, if $H$ is a $p$-group then
$\ed(G; p) \geq \rank(H)$. 

\smallskip
(2) More generally, if the identity component $G^0$ of $G$ is semisimple
and the centralizer $C_{G^0}(H \cap G^0)$ is finite then $\ed(G) \geq 
\rank(H)$. Moreover, if $H$ is a $p$-group then
$\ed(G; p) \geq \rank (H)$. 
\end{thm}

\begin{pf} Apply Theorem~\ref{thm5-5} with $X = V$ = generically free
linear representation of $G$ and $x = 0$.
\end{pf}

\begin{remark} \label{rem.ed=r} Let $G$ be a semisimple algebraic group
and $H = (\bbZ/p^{i_1}\bbZ) \times \dots \times (\bbZ/p^{i_r}\bbZ)$
be an abelian $p$-subgroup of $G$ of rank $r$ satisfying the assumptions
of Theorem~\ref{thm5-6}(1). Then $\ed(G; p) \geq r$ and, in particular,
$\ed(V; p) \geq r$ for any generically free linear representation of $G$; cf. 
Lemma~\ref{lem.ed-p}(b).

Moreover, there exists an irreducible $G$-variety $X$ such that
$\ed(X; p) = r$. Indeed, let $W = \bbA^r$ be a faithful
representation of $H$, where the $i$th cyclic factor of $H$ acts
by a faithful character on the $i$th coordinate of $\bbA^r$, and trivially on
all other coordinates.  Let $X = G \times_H W$ be 
the induced $G$-variety. Since
$X$ is the quotient of the smooth variety $G \times W$ by the free $H$-action
$h(g, w) = (gh^{-1}, hw)$, $X$ is smooth and $\dim X=\dim G+r$.
By our construction the point $x = (1_G, 0_W)$ is fixed by $H$.
Theorem~\ref{thm5-5}(1) shows
that $\ed(X; p) \geq r$; on the other hand, $\dim X-\dim G=r$, and hence,
$\ed(X; p)=r$.

The same construction goes through if $G$ and $H$ satisfy the assumptions
of Theorem~\ref{thm5-6}(2), except that in this case $X$ will be primitive
and not necessarily irreducible.
\end{remark}

%%%%%%%%%%%%%%%%%%%%%%%%%%%

\section{Applications}
\label{sect8}

We now want to apply Theorem~\ref{thm5-6} to specific groups $G$.
In most cases we will always choose $H$ to be an elementary abelian 
$p$-subgroup of $G$.
Note that the theorem does not apply if $H$ is contained in a subtorus $T$
of $G$ because in this case the centralizer of $H$ contains $T$ and, hence, is
infinite. Thus we are interested in nontoral elementary abelian 
$p$-subgroups of $G$. These subgroups have been extensively studied; 
see, e.g., \cite{adams}, \cite{borel}, \cite{bs},
\cite{cs}, \cite{griess}, \cite{wood}.

Before we proceed with the applications, we make two additional remarks.
First of all, for the purpose of applying Theorem~\ref{thm5-6} 
we may restrict our attention to maximal elementary abelian subgroups 
of $G$. Indeed, we lose nothing if we replace $H$ by a larger (with respect 
to containment) elementary abelian subgroup; this will only have
the effect of making the centralizer smaller and improving
the resulting bound on $\ed(G)$.
Secondly, a nontoral elementary abelian subgroup of $G$, even 
a maximal one, may have an infinite centralizer and, hence, not be suitable
for our purposes. Thus our task is to find maximal elementary abelian 
subgroups of $G$ with finite centralizers.

We shall assume that $k$ is an algebraically closed field throughout
this section; cf. Remark~\ref{rem6.05}.

\subsection*{Orthogonal groups}

\begin{thm} \label{thm6.1-1a}
\begin{enumerate}
\item $\ed(O_n; 2) \geq n$ for every $n \geq 1$.
\item $\ed(SO_n; 2) \geq n-1$ for every $n \geq 3$.
\item $\ed(PO_n; 2) \geq n-1$ for every $n \geq 3$.
\end{enumerate}
\end{thm}

\begin{pf} 
Apply Theorem~\ref{thm5-6} with 

\smallskip
(1) $H \simeq (\bbZ/2\bbZ)^n$ = the diagonal subgroup of $G = O_n$.

\smallskip
(2) $H \simeq (\bbZ/2\bbZ)^{n-1}$ = the diagonal subgroup of $G = SO_n$.

\smallskip
(3) $H \simeq (\bbZ/2\bbZ)^{n-1}$ = the diagonal subgroup of $G = PO_n$.
\end{pf}

\begin{remark} For alternative proofs of (1) see \cite[Theorem 10.3 and
Example 12.6]{r2}. For alternative proofs of (2) see \cite[Theorem 10.4 and
Example 12.7]{r2}. (Note that equality holds in both cases.)
The inequality (3) is new to us. 
\end{remark}

\subsection*{Projective linear groups}

The essential dimension of $\PGL_n$ is closely related to the structure of
central simple algebras of degree $n$; we begin by briefly recalling this
connection. 

We shall say that a field extension $K/F$ is {\em prime-to-}$p$ if
it is a finite extension of degree prime to $p$.

\begin{defn} \label{def.tau}
(a) Let $F$ be a field and let $A$ be a finite-dimensional $F$-algebra.
We will say that $A$ is defined over $F_0$ if there exists an 
$F_0$-algebra $A_0$ such that $A \simeq A_0 \otimes _{F_0} F$ (as 
$F$-algebras). Equivalently, $A$ is defined over $F_0$ if there exists
an $F$-basis $e_1, \ldots, e_d$ of $A$ such that $$e_i e_j = \sum_{h=1}^d
c_{ij}^h e_h $$ and every structure constant $c_{ij}^h$ is
contained in $F_0$.

\smallskip
(b) $\tau(A)$ is defined as the minimal value of $\trdeg_k(F_0)$. Here 
the minimum is taken over all subfields $F_0$ of $F$ such that $k \subset F_0$
and $A$ is defined over $F_0$.

\smallskip
(c) Let $p$ be a prime. Then $\tau(A; p)$ is defined as the minimal 
value of $\tau(A \otimes_F K)$, where $K$ ranges over prime-to-$p$ extensions
of $F$.
\end{defn}

\begin{example} If $A = \Mn(F)$ then $\tau(A) = 0$, since $A = \Mn(k) 
\otimes_k F$. 
\end{example}

\begin{lem} \label{lem.pgln}
\begin{enumerate}
\item $\ed(\PGLn)$ is the maximal value of  $\tau(A)$ as $A$ ranges over all
central simple 
algebras of degree $n$ containing $k$ as a central subfield.

\item $\ed(\PGLn)$ is the maximal value of  $\tau(D)$ as $D$ ranges over all
division algebras of degree $n$ containing $k$ as a central subfield.

\item $\ed(\PGLn; p)$ is the maximal value of  $\tau(A; p)$ 
as $A$ ranges over all central simple
algebras of degree $n$ containing $k$ as a central subfield.

\item $\ed(\PGLn; p)$ is the maximal value of  $\tau(D; p)$ 
as $D$ ranges over all
division algebras of degree $n$ containing $k$ as a central subfield.

\item $\ed(\PGLn; p) = \ed(\PGL_{p^r}; p)$, where $p^r$ is 
the highest power of $p$ dividing $n$.

\item $\ed(\PGLn; p) = 0$ if $n$ is not divisible by $p$.

\item $\ed(\PGL_p; p) = 2$.
\end{enumerate}
\end{lem}

\begin{pf} (1) and (2) are proved in \cite[Lemma 9.2]{r2}. (3) and (4)
follow from \cite[Lemma~9.1, Proposition~8.6 and Theorem~8.8(a)]{r2}.

\smallskip
(5) Suppose $n = p^r m$, where $m$ is not divisible by $p$. 
If $D$ is a division
algebra of degree $p^r$ with center $F$ and $A = \Mat_m(D)$
then $\tau(A) = \tau(D)$; see
\cite[Lemma 9.7]{r2}. Thus for any prime-to-$p$ extension $K/F$, 
we have $\tau(A \otimes_F K) = \tau(D \otimes_F K)$. By part
(3) the maximal value of the left hand side (over all $D$ and $K$)
is $\leq \ed(\PGL_n; p)$. On the other hand, by part (4), the maximal
value of the right hand side is $\ed(\PGL_{p^r}; p)$. Thus 
$\ed(\PGL_{p^r}; p) \leq \ed(\PGLn; p)$.

Conversely, given any division algebra $D$ of degree
$n$ with center $F$, there exists a prime-to-$p$ extension $K/F$
such that $D \otimes _F K = \Mat_m(D_0)$, where $D_0$ is a division
algebra of degree $p^r$ with center $K$; see \cite[Theorem 3.1.21]{rowen}.
Thus by part (4)
\[ \tau(D; p) \leq 
\tau(D \otimes_F K; p) = \tau(D_0 ; p) \leq \ed(\PGL_{p^r}; p) \; .\]
Taking the maximum over all $D$ and using part (4) once again, 
we obtain $\ed(\PGLn; p) \leq 
\ed(\PGL_{p^r}; p)$, as desired.

\smallskip
(6) Follows from part (5) with $r = 0$.

\smallskip
(7) It is enough to show $\tau(D; p) = 2$ for every division algebra $D$
of degree $p$. To show $\tau(D; p) \geq 2$, note that  
for any prime-to-$p$ extension $K/F$, $D \otimes _F K$ is 
a division algebra; see~\cite[Corollary 3.1.19]{rowen}. By Tsen's theorem, 
$\tau(D \otimes_F K) \geq 2$; see \cite[Lemma 9.4(a)]{r2}. 
This proves $\tau(D; p) \geq 2$.

On the other hand,
by a theorem of Albert, there exists a prime-to-$p$ extension $K/F$ 
such that $D' = D \otimes_F K$ is a cyclic division algebra. 
Then by \cite[Lemma 9.4(b)]{r2} $\tau(D') \leq 2$ and hence,
$\tau(D; p) \leq 2$. 
\end{pf}

The following inequality is a consequence of \cite[Theorem 16.1(b)]{joubert}
and Lemma~\ref{lem.pgln}(3) above.  

\begin{thm} \label{thm.pgln}
$\ed(\PGL_{p^r}; p) \geq 2r$.
\end{thm}

We will now give an alternative proof based on Theorem~\ref{thm5-6}. 
In fact, we will prove a slightly stronger result; see Theorem~\ref{thm.sl}. 
We begin with the following elementary construction.

\begin{defn} \label{def.phi_A} 
Let $A$ is an abelian group of order $n$. 
and let $V = k[A]$ be the group algebra of $A$.  

\smallskip
(a) The regular representation $P \colon A \lra \GL(V) = \GL_n$
is given by $a \mapsto P_a \in \GL(V) = \GL_n$,
where
\[ P_a\Bigl(\sum_{b \in A} c_b b\Bigr) = \sum_{b \in A} c_b ab  \]
for any $a \in A$ and $c_b \in k$. 
  
\smallskip
(b) The representation $D \colon A^{\ast} \lra \GL(V) = \GL_n$ is
defined by $\chi \mapsto D_{\chi} \in \GL(V)$, where
\[ D_{\chi}\Bigl( \sum_{a \in A} c_a a\Bigr) = \sum_{a \in A} c_a \chi(a) a  \]
for any $\chi \in A^{\ast}$ and $c_a \in k$. 

\smallskip
Note that in the basis $\{ a  \mid a \in A \}$ of $V$,
each $P_{a}$ is represented by a permutation matrix and each $D_{\chi}$ is
represented by a diagonal matrix; this explains our choice of the letters
$P$ and $D$.
\end{defn}

\begin{lem} \label{lem.formulas}
Let $A$ be a finite abelian group, $a, b \in A$ and
$\chi, \mu \in A^{\ast}$. Then

\smallskip
(a) $D_{\chi} P_a = \chi(a) P_a D_{\chi}$. 

\smallskip
(b) $(P_a D_{\chi}) (P_b D_{\mu}) (P_a D_{\chi})^{-1} = \chi(b) \mu^{-1}(a) 
(P_b D_{\mu})$
\end{lem}

\begin{pf} 
Part (a) can be verified directly from Definition~\ref{def.phi_A}.
Part (b) is an immediate consequence of part (a).
\end{pf}

\begin{lem} \label{lem.SL}
Suppose $A$ is an abelian group of order $n$ such that its $2$-Sylow
subgroup is either (i) non-cyclic or (ii) trivial (the latter possibility
happens when $n$ is odd).  Then

\smallskip
(a) $P_a \in \SL_n$ for every $a \in A$ and

\smallskip
(b) $D_\chi \in \SL_n$ for every $\chi \in A^{\ast}$.
\end{lem}

\begin{pf} (a) Recall that $P_a$ is a permutation matrix representing 
the permutation
$\sigma_a \colon A \lra A$ given by $b \lra ab$. Thus $\det(P_a) =
(-1)^{{\rm sign}(\sigma_a)}$, and we only need to show $\sigma_a$ is
even.

Assume, to the contrary, that $\sigma_a$ is odd. 
Let $m$ be the order of $a$.  Since $\sigma_a$ is a product of
$\frac{n}{m}$ disjoint $m$-cycles, both $\frac{n}{m}$ and $m-1$ 
are odd.  In particular, $m$ and, hence,
$n$ is even. Thus assumption (ii) fails. On the other hand, since
$\frac{n}{m} = [A : \lf < a \r>]$ is odd, the Sylow 2-subgroup
of $A$ is contained in $\lf< a \r>$ and, thus assumption (i) fails.
This contradiction proves that $\sigma_a$ is an even permutation.

\smallskip
(b) Suppose $\chi$ is an element of $A^*$ of order $m$ 
and let $\zeta_m$ be a primitive $m$-th root of unity.
The matrix $D_{\chi}$ is diagonal with entries $\chi(a)$, as $a$
ranges over $A$; here $\chi(a)$ assumes the value $(\zeta_m)^i$ exactly
$\frac{n}{m}$ times for each $i = 0, 1, \dots, m-1$.
Hence,
\[
\det(D_{\chi}) =
\Bigl(\prod_{i=0}^{m-1}(\zeta_m)^i\Bigr)^{\frac{n}{m}} =
{(\zeta_m)\vphantom{\Big)}}^{m\cdot\frac{m-1}{2}\cdot\frac{n}{m}}\ .
\]
Assume, to the contrary that $\det(D_{\chi})\neq 1$. Then
both $\frac{n}{m}$ and $m-1$
are odd.  Arguing as in part (a), we conclude that $n = |A^{\ast}|$ is
even and the Sylow 2-subgroup of $A^{\ast}$ is cyclic.
Since $A$ and $A^*$ are isomorphic, this contradicts our
assumption. Hence, $\det(D_{\chi}) =1$, as claimed.
\end{pf}

\begin{defn} \label{def.phi}
Assume $A$ is an abelian group of order $n$, $e$ is an integer dividing
$n$ and $\zeta_e$ is a primitive $e$th root of unity.

\smallskip
(i) Let $\phi_n \colon A \times A^{\ast} \lra \PGLn$ be the map
of sets given by $\phi_n(a, \chi)$ = image of $P_a D_{\chi}$ in $\PGLn$.
We define $H_n$ as the image of $\phi_n$ in $\PGLn$.

\smallskip
(ii) Suppose $A$ satisfies the conditions of Lemma~\ref{lem.SL}.
Then we define $\phi_e \colon A \times A^{\ast} \lra 
\SL_n/ \lf< \zeta_e I_n \r> $ by the formula
$\phi_e(a, \chi) = P_a D_{\chi}$ (mod $< \zeta_e I_n \r>$).
We define $H_e$ as the image of $\phi_e$ in
$\SL_n/ \lf< \zeta_e I_n \r>$. 

Note that $\SL_n/\lf<\zeta_n I_n \r> = \PGLn$. If
$A$ satisfies the conditions of Lemma~\ref{lem.SL}
then  the two definitions of $\phi_n$ (and thus $H_n$) coincide.
\end{defn}

\begin{lem} \label{lem.PD} 
In the assumptions of Definition~\ref{def.phi},

\smallskip
(i) $H_n$ is a subgroup of $\PGLn$ and
$\phi_n$ is an isomorphism between $A \times A^{\ast}$ and $H_n$;

\smallskip
(ii) $H_e$ is a subgroup of $\SL_n/ \lf< \zeta_e I_n \r>$ and
$\phi_e$ is an isomorphism between $A \times A^{\ast}$ and $H_e$,
provided that $A$ satisfies the conditions of Lemma~\ref{lem.SL}
and the exponent of $A$ divides $e$.
\end{lem}

The lemma says, in particular, that, if $e$ is divisible 
by the exponent of $A$ then
$H_e$ is a subgroup of $\SL_n/ \lf< \zeta_e I_n \r>$ whenever $H_e$ is 
defined.  (Note that in part (i),
$e=n = |A|$ is necessarily divisible by the exponent of $A$.)

\begin{pf} By Lemma~\ref{lem.formulas}(a),
$P_a$ and $D_{\chi}$ commute modulo $\lf< \zeta_{n}I_n \r>$ in case
(i) and modulo $\lf< \zeta_{e}I_n \r>$ in case (ii).
The lemma is an easy consequence of this fact.
\end{pf}

In the sequel we shall assume that $A$ is an abelian $p$-group of order
$p^r$ and $e=p^i$, where $1\le i\le r$ is chosen so that $e$ is divisible
by the exponent of $A$. 
Note that under these assumptions $H_e$ is always well-defined and is
a subgroup of $SL_n/\lf< \zeta_e I_n \r>$. (Indeed,
if the conditions of Lemma~\ref{lem.SL} fail to be satisfied then
$p=2$, $A$ is cyclic and hence, $e = n$, so that $H_e$ is
given by Definition~\ref{def.phi_A}(i).)

\begin{lem} \label{lem.self-c}
Let $A$ be an abelian $p$-group of order $n = p^r$, let $e=p^i$ with
$1\le i \leq r$. Assume the exponent of $A$ divides $e$.
Let $\pi \colon SL_n/ \lf< \zeta_e I_n \r> \lra \PGLn$ 
be the natural projection, let $H = \pi^{-1}(H_n)$ and 
let $K = \Ker(\pi)$ be the center of $SL_n/ \lf< \zeta_e I_n \r>$. Then

\smallskip
(a) $H = H_e  \times K \simeq A \times A^{\ast} \times (\bbZ/p^{r-i}\bbZ)$.

\smallskip
(b) $H_n$ is self-centralizing in $\PGLn$,

\smallskip
(c) $H$ is self-centralizing in $SL_n/\lf< \zeta_e I_n \r>$.
\end{lem}

\begin{pf} (a) The surjective homomorphism $\pi|_{H} \colon H \lra H_n$ 
splits: the complement of $K$ in $H$ is $H_e$.
Since $K$ is central, part (a) follows.

\smallskip
(b) Denote the centralizer of $H_n$ in $\PGLn$ by $C(H_n)$.
Lemma~\ref{lem.formulas}(b) shows that
for every $b \in A$ and $\mu \in A^{\ast}$, 
the matrix $P_b D_{\mu}$ spans a one-dimensional representation space for the
conjugation action of $H_n$ on $\Mn(k)$; moreover, 
$H_n$ acts on these $| H_n |$ spaces by distinct characters. 
Since there are $n^2 = p^{2r}$ of these spaces and 
$\dim(\Mn) = | H_n | = n^2$, we conclude that $\Mn(k)$ decomposes as a direct 
sum of these one-dimensional representations.
Any $g \in C(H_n) \subset \PGLn$ is represented by a non-zero matrix lying 
in one of them, i.e., by a non-zero constant multiple of
$P_a D_{\chi}$ for some $a \in A$ and $\chi \in A^{\ast}$. This shows
that $C(H_n) = H_n$ in $\PGLn$, as claimed.

\smallskip
(c) Denote the centralizer of $H$ in $SL_n/\lf< \zeta_e I_n \r>$ by $C(H)$.
Since $H$ is abelian, $H \subset C(H)$. On the other hand, in view of part
(b), $C(H) \subset \pi^{-1}(C(H_n)) = \pi^{-1}(H_n) = H$. 
\end{pf}

\begin{thm} \label{thm.sl}
\[
\ed(\SL_{p^r}/\lf< \zeta_{p^i} I_{p^r} \r>;p) \geq
\begin{cases}
2r+1 & \text{if $i=1,\dots,r-1$}\ ,\\
2r & \text{if $i=r$} \ .
\end{cases}
\]
\end{thm}

Note that if $i=r$ then
$\SL_{p^r}/ \lf< \zeta_{p^i} I_{p^r} \r>  = \PGL_{p^r}$, and 
we obtain the bound of Theorem~\ref{thm.pgln}.

\begin{pf}
Applying Lemma~\ref{lem.self-c} to $A=(\bbZ/p\bbZ)^r$; we obtain 
a finite abelian self-centralizing $p$-subgroup 
$H \subset\SL_{p^r}/ \lf< \zeta_{p^i} I_{p^r} \r>$.
By part (a) $\rank(H) = 2r+1$
if $1\le i < r$ and $2r$ if $i=r$. The desired inequalities now follow from
Theorem~\ref{thm5-6}. 
\end{pf}

\begin{remark} \label{rem.pgln} 
One can show that any abelian $p$-subgroup of $\PGL_{p^r}$ with a finite
centralizer has rank $\leq 2r$. Thus the lower bounds of Theorem~\ref{thm.sl}
cannot be improved by this method.
\end{remark}

\subsection*{Spin groups}
We will now apply Theorem~\ref{thm5-6} to obtain lower bounds
on the essential dimension of some spin groups. Elementary abelian subgroups
of $\Spin_n$ are described in some detail in \cite{wood}. In particular, 
if $p$ is an odd prime then every elementary abelian $p$-group is toral 
(see~\cite[Section 2.2]{serrepp}, \cite[Theorem 5.6]{wood}, 
or~\cite[(2.22)]{griess}) and thus is not suitable
for our purposes. We shall therefore concentrate on elementary
abelian 2-subgroups.

Recall that $\Spin_n$ fits into an exact sequence 
$$ \{1\} \lra \{ -1, 1 \} \lra \Spin_n \stackrel{f}{\lra} SO_n \lra \{1 \}
\; , $$ 
where $\{-1, 1\}$ is the central subgroup of $\Spin_n$.  Let 
$D \simeq (\bbZ/2\bbZ)^{n-1}$ be the diagonal subgroup of $SO_n$
and let $D' = f^{-1}(D) \subset \Spin_n$. We want to construct elementary
abelian 2-subgroups of $D'$. (Note that every
elementary abelian 2-subgroups of $\Spin_n$ is conjugate to a subgroup of
$D'$; see~\cite[Theorem 5.6]{wood}.) 

Recall that a doubly even code $L$ of length $n$ is a vector subspace of
$(\bbZ/2\bbZ)^n$ with the property that the weight of every element of $L$ 
is divisible by 4.
(Here the weight of an element of  $(\bbZ/2\bbZ)^n$ is defined as the number 
of 1s among its coordinates.) We shall say that an $m \times n$-matrix
over $\bbZ/2\bbZ$ is a generator matrix for $L$ if its rows span $L$
as a $\bbZ/2\bbZ$-vector space.

Doubly even codes of length $n$ are in 1---1
correspondence with elementary abelian 2-subgroups of $D'$ containing $-1$;
this is explained in \cite[Sections 1 and 2]{wood}; see 
also~\cite[Section 7]{steinberg}.
Explicitly, let $E_n$ be the (index 2) subgroup of $(\bbZ/2\bbZ)^n)$ 
consisting of all codewords of even weight.  Consider 
the group isomorphism $\phi\colon (E_n, +) \lra (D, \cdot)$ given by 
\begin{equation} \label{e6.11}
\phi(i_1, \ldots, i_n) =
\left( \begin{matrix} (-1)^{i_1} & 0 & \hdots & 0 \\
                               0 & (-1)^{i_2} & \hdots & 0 \\
                               \hdots & \hdots & \hdots & \hdots \\
                               0 & 0 & \hdots & (-1)^{i_n} \end{matrix}
\right) \end{equation}
If $L \in E_n$ is a doubly even code of dimension $d$ then
$\phi(L)$ is an elementary abelian 2-subgroup of $SO_n$ of rank $d$. The 
preimage $H = f^{-1}(\phi(L))$ of this subgroup in $\Spin_n$
is thus an elementary abelian 2-subgroup of rank $d+1$. Note that
by~\cite[Theorem 2.1]{wood}, every elementary abelian
subgroup of $D'$ containing $-1$ is obtained in this way.

Recall that not every elementary abelian 2-subgroup is good for our purposes;
in order to apply Theorem~\ref{thm5-6}, we need to construct one whose
centralizer is finite. Clearly the group $H= f^{-1}(\phi(L))$ has 
a finite centralizer in $\Spin_n$ if and only if its image $f(H) = \phi(D)$ 
has a finite centralizer in $SO_n$.

\begin{lem} \label{lem7.11} 
Let $L$ be a doubly even code of length $n$ and let
\[ \phi\colon E_n \lra SO_n \] be as in \eqref{e6.11}. 
Then $\phi(L)$
has a finite centralizer in $SO_n$ if and only if a generator matrix of $L$ 
has distinct columns. 
\end{lem} 

\begin{pf} The map $\phi|_L$ may be viewed as an orthogonal 
representation of $L \simeq (\bbZ/2\bbZ)^d$. This representation is given
to us as a direct sum of characters $\chi_1, \dots, \chi_n\colon L \lra 
\bbG_\ml$,
where $\chi_j(i_1, \ldots, i_n) = (-1)^{i_j}$. 
Note that a generator matrix of $L$ has
distinct columns if and only if these characters are distinct.
If the characters are distinct then by Schur's Lemma
the centralizer of $\phi(L)$ in
$SO_n$ consists of diagonal matrices and, hence, is finite. On the other
hand, if two of these characters are equal
then the centralizer of $\phi(L)$ contains a copy of $SO_2$ and, hence, is
infinite.
\end{pf}

We are now ready to state our main result on spin groups.

\begin{thm} \label{thm7.11} 
$\ed(\Spin_n; 2) \geq [\dfrac{n}{2}]+ 1$ for every 
$n \equiv 0$, $1$ or $-1 \pmod 8$. 
\end{thm}

\begin{pf} The above discussion shows that it is sufficient to
construct a doubly even code $L$ of length $n$ and dimension $[n/2]$ all
of whose columns are distinct. 

We now exhibit such codes in the three cases
covered by the theorem. 
Let $0_i$ (respectively, $J_i$) denote, the $i$-tuple 
of zeros (respectively, the $i$-tuple of ones) in $(\bbZ/2\bbZ)^i$.
One can now check directly that each of the following codes is doubly even
of dimension $[\dfrac{n}{2}]$; moreover, in each case the generator matrix 
(for the generating set given below) has distinct columns.

\smallskip
$n = 8m$. $L = \lf<(a, a), (0_{4m}, J_{4m})\r>$, where $a$ ranges over
all elements of $(\bbZ/2\bbZ)^{4m}$ of even weight.

\smallskip
$n = 8m+1$. $L = \lf<(0_1, a, a), (0_{4m+1}, J_{4m})\r>$, where $a$ ranges over
all elements of $(\bbZ/2\bbZ)^{4m}$ of even weight.

\smallskip
$n = 8m-1$. $L = \lf<(a, a, 0_1), (0_{4m-1}, J_{4m})\r>$, where $a$ ranges over
all elements of $(\bbZ/2\bbZ)^{4m-1}$ of even weight.

\smallskip
This completes the proof of the theorem.
\end{pf}

\begin{remark} \label{rem.spin2-6} 
Recall the following exceptional isomorphisms of classical algebraic groups:
\[ \begin{array}{l} 
\Spin_2 \simeq (\bbG_\ml)^2 
\, , \\
\Spin_3 \simeq \SL_2 \, , \\
\Spin_4 \simeq \SL_2 \times \SL_2 \, , \\
\Spin_5 \simeq {\rm Sp}_4\, , \; \text{and} \\
\Spin_6 \simeq \SL_4 \, .
\end{array} \]
(This phenomenon is caused by
the fact that while the Dynkin diagrams of 
types $A_n$, $B_n$, $C_n$, and $D_n$ 
are distinct for large $n$, for small $n$ there are some overlaps.) 
We conclude that all of these groups are special (see
\cite[Section 5]{grothendieck}, \cite[Section 2.6]{pv}) and thus
\[ \ed(\Spin_n) = 0 \; \, 
\text{for every} \; \, 2 \leq n \leq 6  \; ; \]
(see \cite[Section 5.2]{r2}). This shows 
that the condition $n \equiv 0$, $1$ or $-1 \pmod 8$ is not 
as arbitrary as it may seem at first glance.
\end{remark}

\begin{remark} \label{rem.rost}
The following results are due to M. Rost~\cite{rostspin}:
\[ \begin{array}{l}
\ed( \Spin_7) = 4 \\
\ed( \Spin_8) = 5 \\
\ed( \Spin_9) = 5 \\
\ed( \Spin_{10}) = 4 \\
\ed( \Spin_{11}) = 5 \\
\ed( \Spin_{12}) = 6 \\
\ed( \Spin_{13}) = 6 \\
\ed( \Spin_{14}) = 7 \; . \end{array} \]
The proofs rely on the properties of quadratic forms 
of dimension $\leq 14$. In particular,
our bound is sharp for $n = 7$, $8$ and $9$.  On a lighter 
note, our bound is also sharp for $n = 1$, since $\Spin_1 = \bbZ/2\bbZ$
and $\ed(\bbZ/2\bbZ) = 1$. 
\end{remark}

\subsection*{Exceptional groups}

\begin{thm} \label{thm6.1-1}
\begin{enumerate}
\item $\ed(G_2; 2) \geq 3$.
\item $\ed(F_4; 2) \geq 5$.
\item $\ed(F_4; 3) \geq 3$.
\item $\ed(3E_6; 3) \geq 4$. Here $3E_6$ denotes the simply connected 
group of type $E_6$ over $k$.
\item $\ed(2E_7; 2) \geq 7$. Here $2E_7$ denotes the simply connected 
group of type $E_7$ over $k$.
\item $\ed(E_7; 2) \geq 8$.  Here $E_7$ denotes the adjoint $E_7$.
\item $\ed(E_8; 2) \geq 9$.
\item $\ed(E_8; 3) \geq 5$.
\item $\ed(E_8; 5) \geq 3$.
\end{enumerate}
\end{thm}

\begin{pf} In each case we exhibit an abelian subgroup $H$ 
with a finite centralizer, then appeal to Theorem~\ref{thm5-6}. 

\smallskip
(1) Let $\Oct$ be the split octonion algebra generated by $i$, $j$, and $l$,
as in \cite[pp.~16--17]{jacobson}. We can identify
$G_2 \subset \GL_8$ with the automorphism group of $\Oct$. Now
let $H = \lf<\alpha, \beta, \gamma\r> \simeq (\bbZ/2\bbZ)^3$, where
$$ \begin{array}{lll} \alpha(i) = -i \,  & \alpha (j) = j  & \alpha(l) = l \\
 \beta(i) = i \,  & \beta (j) = -j  & \beta(l) = l \\
 \gamma(i) = i \,  & \gamma (j) = j  & \gamma(l) = -l \; . \end{array} $$
To prove that $H$ is self-centralizing, note that
the representation of $H$ on $\Oct$ (viewed as an 8-dimensional vector 
space) is a direct sum of 8 distinct characters; 
cf.~\cite[Table I, p. 257]{griess} or \cite[p. 252]{cs}.

\smallskip
(2) A self-centralizing $H = (\bbZ/2\bbZ)^5 \subset F_4$ is described 
in~\cite[(7.3)]{griess}.

\smallskip
(3) A self-centralizing $H = (\bbZ/3\bbZ)^3 \subset F_4$ is described 
in~\cite[(7.4)]{griess}.

\smallskip
(4) Use the maximal $H = (\bbZ/3\bbZ)^4$ of $3E_6$ 
described in \cite[(11.13)(i)]{griess}; see also~\cite{cs}.
Note that by \cite[(11.13)(i)]{griess} $H$ has a finite normalizer in
$3E_6$; hence, its centralizer is finite as well.

\smallskip
(5) Let $u$ be an element of order 4 in $2E_7$
whose centralizer $C(u)$ is isomorphic to $\SL_8/(\pm I_8)$; see
\cite[bottom of p. 283]{griess}. (According to the notational 
conventions of \cite[(2.14)]{griess}, $u$ is an element of type
{\bf 4A}.) Note that under the identification
$C(u)\isomo\SL_8/(\pm I_8)$, the element $u$ corresponds to the
central element of order $4$ in $\SL_8/(\pm I_8)$ which is represented
by the identity matrix $I_8$.

By Lemma~\ref{lem.self-c}(a), with $p =2$, $r = 3$, $e= 2$ and 
$A = (\bbZ/2\bbZ)^3$, the group
$C(u) = \SL_8/(\pm I_8)$ contains a self-centralizing finite abelian subgroup 
$H \simeq (\bbZ/2\bbZ)^6 \times (\bbZ/4\bbZ)$, where the
$\bbZ/4\bbZ$-factor is the center of $C(u)$, i.e., is equal to
$\lf< u \r>$. Moreover, $H$ is self-centralizing 
in $C(u)$. Since $u \in H$, we conclude that $H$ is self-centralizing in
$2E_7$. Applying Theorem~\ref{thm5-6} to $H$, we obtain
the desired inequality $\ed(2E_7) \geq \rank(H) = 7$.

\smallskip
(6) A self-centralizing subgroup $H = (\bbZ/2\bbZ)^8$ of $E_7$ 
is described in \cite[Theorem~9.8(ii)]{griess}; see also~\cite{cs}.

\smallskip
(7) $E_8$ has a maximal elementary abelian subgroup $H \simeq (\bbZ/2\bbZ)^9$ 
called a ``type 1 subgroup"; see~\cite{adams}, \cite[(2.17)]{griess} 
and \cite{cs}. By \cite[(2.17)]{griess} this subgroup has a finite normalizer.
Hence, its centralizer is finite as well. (In fact, one can
show that $H$ is self-centralizing; see~\cite[p. 258]{griess}).

\smallskip
(8)--(9) $E_8$ contains self-centralizing subgroups 
$H_1 \simeq (\bbZ/3\bbZ)^5$; 
and $H_2 = (\bbZ/5\bbZ)^3$; see~\cite[(11.5) and (10.3)]{griess} 
\end{pf}

\begin{remark} \label{rem6.1-1}
Alternative proofs of the inequalities (1), (2) and (3)
can be found in \cite[12.14 and 12.15]{r2}. In fact, equality holds 
in all three cases: in the case of (1) this is proved 
in \cite{r2}, for (2) and (3) this was observed by 
J.-P. Serre~\cite{serrepc1}. Moreover,
V. E. Kordonsky~\cite{kordonsky} has shown that
$\ed(F_4) \leq 5$ (and thus $\ed(F_4) = 5$).

One can show, by modifying the proof of
\cite[Proposition 11.7]{r2} (or, alternatively, 
of~\cite[Theorem 9]{kordonsky}) that $\ed(3E_6; 3) \leq \ed(F_4; 3) + 1 = 4$,
so that inequality (4) is sharp as well.
We do not know the exact value of $\ed(3E_6)$; however, Kordonsky has 
shown that $\ed(3E_6) \leq 6$; see~\cite[Section 4.2]{kordonsky}.  
Thus $\ed(3E_6) = 4$, $5$ or $6$.
We remark that alternative proofs of (4) were recently 
shown to us by M. Rost 
and by R. S. Garibaldi~\cite{garibaldi}. 

An alternative proof of part (9) is based on Lemma~\ref{lem.ci} and
the existence of a nontrivial Rost invariant
$H^1(\, \cdot \,, E_8) \lra H^5(\, \cdot \,, \bbZ/5\bbZ)$;
see \cite[7.3]{serrepp} or \cite[(31.40) and (31.47)]{boi}. 
M. Rost has pointed out to us that, in fact, $\ed(E_8; 5) = 3$.

We do not know whether or not inequalities (5)--(8) are sharp.
Regarding (5), we remark that by a theorem of Kordonsky $\ed(2E_7) \leq 9$ 
(see \cite[Theorem 10]{kordonsky}); thus $\ed(2E_7)$ and $\ed(2E_7; 2)$ 
are equal to $7$, $8$ or $9$.

To the best of our knowledge, the inequalities (5)--(8) are new. 
\end{remark}

\subsection*{A wish list for cohomological invariants}

\begin{remark} \label{rem.ci2} 
Some of the lower bounds of this section allow alternative proofs based
on the existence of certain cohomological invariants; see~Lemma~\ref{lem.ci}. 
For example, Theorem~\ref{thm6.1-1a}(1) follows from the existence
of a non-trivial cohomological invariant $H^1(\, \cdot \,, O_n)
\lra H^n(\, \cdot \,, \bbZ/2\bbZ)$ (namely, the $n$th Stiefel---Whitney 
class, see \cite[Section 6.3]{serrepp}),
Theorem~\ref{thm6.1-1}(2) follows from the existence
of the cohomological invariant of $H^1(\, \cdot \,, F_4) \lra
H^5(\, \cdot \,, \bbZ/2\bbZ)$ (see~\cite[Section 9.2]{serrepp}), 
Theorem~\ref{thm6.1-1}(3) follows from
the existence of the Serre---Rost invariant  $H^1(\, \cdot \,, F_4) \lra
H^5(\, \cdot \,, \bbZ/3\bbZ)$ (see~\cite[Section 9.3]{serrepp}), etc.

Other inequalities cannot be proved in this way because the needed
cohomological invariants are not known to exist. On the other 
hand, these bounds suggest that there may exist cohomological invariants 
of the types listed below.
(Here by a mod $p$ invariant of $G$-varieties in $H^d$ 
we shall mean a cohomological invariant 
$H^1(\, \cdot \,, G) \lra H^d (\, \cdot \, , M)$ in the sense of 
Definition~\ref{def.ci}, with $M$ $p$-torsion.) 

\begin{enumerate}
\item (cf.\ Theorem~\ref{thm.pgln})
A mod $p$ invariant of $\PGL_{p^r}$-varieties in $H^{2r}$.

In the case $p=r=2$ an invariant of this type was recently constructed
by J.-P. Serre~\cite{serrepc2} (see also \cite{rostprep}). 

\item (cf.\ Theorem~\ref{thm7.11})
A mod $2$ invariant of $\Spin_n$-varieties in $H^{[n/2]+1}$ for
$n \equiv  0, \pm 1 \pmod 8$.

For $n = 7$, $8$ and $9$ such invariants were recently constructed
by M.~Rost~\cite{rostspin}. 

\item (cf.\ Theorem~\ref{thm6.1-1}(6))
A mod $2$ invariant of $E_7$-varieties in $H^8$.

\item (cf.\ Theorem~\ref{thm6.1-1}(7))
A mod $2$ invariant of $E_8$-varieties in $H^9$.

\item (cf.\ Theorem~\ref{thm6.1-1}(8))
A mod $3$ invariant of $E_8$-varieties in $H^5$.
\end{enumerate}

\smallskip
The above-mentioned constructions of Serre and Rost 
represent the only currently known invariants of types 1--5.
\end{remark}
%%%%%%%%%%%%%%%%%%%%%%%%%%%

\section{Simplifying polynomials by Tschirnhaus transformations}
\label{sect9}

Let $E/F$ be a field extension of degree $n$ such that $k \subset F$.
Suppose $E = F(z)$ and 
$$ f_z(t) = t^n + \alpha_1(z)t^{n-1} + \ldots + \alpha_n(z) $$
is the minimal polynomial of $z$ over $F$. 
We are interested in choosing the generator $z$ whose minimal polynomial
has the simplest possible form. More precisely, we want 
$\trdeg_k \, k(\alpha_1(z), \ldots, \alpha_n(z))$ to be as small 
as possible. We shall denote the minimal value of 
$\trdeg_k \, k(\alpha_1(z), \ldots, \alpha_n(z))$ by $\tau(E/F)$.  Note that
$\tau(E/F)$ is the same as $\tau(E)$ given by Definition~\ref{def.tau},
where $E$ is viewed as an $n$-dimensional $F$-algebra.  (We remark 
that $\tau(E/F)$ was denoted by $\ed(E/F)$ in \cite{br1} and \cite{br2}.)

As we explained in the Introduction, a choice of a generator
$z$ (or, equivalently, an isomorphism of fields $E \simeq F[t]/(f_z)$)
is called a Tschirnhaus transformation without auxiliary radicals.
If $E/F$ is given as the root field of a polynomial $f(x) \in F[x]$, i.e.,
$E = F[x]/(f(x))$, then the polynomial $f_z(t)$ is said to be obtained
from $f(t)$ via the Tschirnhaus substitution $x \longmapsto z$. In this setting
we are interested in simplifying the given polynomial $f(t) = f_x(t)$ 
by a Tschirnhaus substitution, where the ``complexity" of a polynomial 
is measured by the
number of algebraically independent coefficients (over $k$). The number 
$\tau(E/F)$ tells us to what extent $f(x)$ can be simplified.

A case of special interest is the generic field extension $L/K$ of degree $n$.
More precisely, $K = k(a_1, \ldots, a_n)$, $L = K[x]/g(x)$, and
\[ f(x) = x^n + a_1 x^{n-1} + \ldots + a_n \; , \]
where $a_1, \ldots, a_n$ are algebraically independent variables over $k$.
The following results are proved in \cite{br1} (see also~\cite{br2}):
$\tau(L/K) = \ed(\Sym_n) \geq [n/2]$ and
$\tau(L/K) \geq \tau(E/F)$, where $E/F$ is any field extension 
of degree $n$.  

The object of this section is to prove Theorem~\ref{thm1.3} stated in the
Introduction. Using the terminology we introduced above,
Theorem~\ref{thm1.3} can be rephrased as follows.

\begin{thm} \label{thm7.1}
Suppose $\dfrac{n}{2} \leq m \leq n-1$, where $m$ and $n$ are positive
integers. 
Let $a_m, \dots, a_n$ be algebraically independent variables over $k$,
$F = k(a_m, \ldots, a_n)$ and $E = F[x]/f(x)$, where 
$$f(x) = x^n + a_m x^{n-m} + \dots + a_{n-1} x + a_n \; . $$
Then $\tau(E/F) = n-m$. 
\end{thm}

Note that $f(x)$ is an irreducible polynomial over $F$ so that $E$ is,
in fact, a field. Indeed, by Gauss' Lemma (see \cite[V.6]{lang}) 
it is enough to check irreducibility over the ring 
$k[a_m, \dots, a_n]$; now we can set $a_m = \dots = a_{n-1}=0$ and 
apply the Eisenstein criterion (see \cite[V.7]{lang}).
Alternatively, the irreducibility of $f(x)$ follows from Lemma~\ref{lem7.2c}
below.

\subsection*{The variety $X_{m,n}$}

Before we can proceed with the proof of Theorem~\ref{thm7.1}, we need
to establish several elementary properties of the variety 
$X_{m, n} \subset \bbA^n$ given by
\begin{equation} \label{e7.1}
X_{m, n} = \{ x = (x_1, \ldots, x_n) \mid s_1(x) = s_2(x) = \dots = 
s_{m-1}(x) = 0 \} \; , 
\end{equation}
where $s_i(x)$ is the $i$th elementary symmetric polynomial 
in $x_1, \ldots, x_n$. Note that $X_{m, n}$ can also be described as
\begin{equation} \label{e7.1-1}
X_{m, n} = \{ x = (x_1, \ldots, x_n) \mid p_1(x) = p_2(x) = \dots = 
p_{m-1}(x) = 0 \} \; , 
\end{equation}
where $p_i(x) = x_1^i + \ldots + x_n^i = 0$; the equivalence 
of the two definitions follows from Newton's formulas.
(Recall that $\Char(k) = 0$ throughout this paper.) 
Note that \eqref{e7.1-1} defines $X_{m, n}$ for every positive integer $m$
(of course, $X_{m, n} = \{ 0 \}$ if $m > n$) and that
the symmetric group $\Sym_n$ acts on $X_{m, n}$ by permuting 
the coordinates $x_1, \ldots, x_n$.

To simplify the exposition, we shall assume that the base field $k$ over which
$X_{m, n}$ is defined, is algebraically closed; we note that
Lemmas~\ref{lem7.2b} and~\ref{lem7.2c} are true without this assumption.

\begin{lem} \label{lem7.2a}
Suppose $x=(x_1, \ldots, x_n) \in X_{m, n}$. Then either 
$x = 0$ or at least $m$ of its coordinates $x_1, \ldots, x_n$ are distinct.
\end{lem}

\begin{pf} It is enough to prove the lemma under the assumption that
$x_i \neq 0$ for every $i = 1, \ldots, n$. Indeed if, say, 
$x_1 = \ldots = x_r =0$ and $x_{r+1}, \ldots, x_n \neq 0$ then 
we can replace $n$ by $n-r$ and $x$ by 
$y= (x_{r+1}, \ldots, x_n) \in X_{i, n-r}$. 

After permuting the coordinates of $x$, we may assume $x_1, \ldots, x_r$ are
distinct and $x_1, \ldots, x_n \in \{x_1, \ldots, x_r\}$. Suppose
$n_1$ of the coordinates $x_1, \ldots, x_n$ are equal to $x_1$,
$n_2$ of them are equal to $x_2$, $\dots$, and $n_r$ of them are equal to
$x_r$. By definition of $X_{m, n}$ we have
$p_1(x) = \ldots = p_{m-1}(x) = 0$ or, equivalently,  
\[ \sum_{i=1}^r n_i x_i^j = 0 \quad \text{for every} \quad j=1,...,m-1 \; .\]
This means that the columns of the Vandermonde matrix
\[ \left( \begin{matrix} x_1 & x_2 & \hdots & x_r \\
                         x_1^2 & x_2^2 & \hdots & x_r^2 \\
                         \dots & \dots & \dots & \dots \\
                         x_1^{m-1} & x_2^{m-1} & \hdots & x_r^{m-1} 
\end{matrix} \right) \]
are linearly dependent. Since we are assuming $x_1, \ldots, x_r$ 
are distinct non-zero elements of $k$, this is only possible 
if $r \geq m$, as claimed. 
\end{pf}

\begin{lem} \label{lem7.2b}
Every non-zero point of $X_{m, n}$ is smooth.
\end{lem}

\begin{pf}
We apply the Jacobian criterion to the system of polynomial equations
$p_1(x) = \dots = p_{m-1}(x) = 0$ defining $X_{m, n}$. The Jacobian matrix
of this system is given by
\[ J(x_1, \ldots, x_n) = \left( \begin{matrix} 1 & 1 & \hdots & 1 \\
                         2x_1 & 2x_2 & \hdots & 2x_n \\
                         3x_1^2 & 3x_2^2 & \hdots & 3x_n^2 \\
                         \dots & \dots & \dots & \dots \\
                  (m-1) x_1^{m-2} & (m-1) x_2^{m-2} & \hdots & (m-1) x_n^{m-2} 
\end{matrix} \right) \; . \]
It is easy to see that this $(m-1) \times n$-matrix
has rank $m-1$ whenever $m-1$ or more of the coordinates
$x_1, \ldots, x_n$ are distinct. By Lemma~\ref{lem7.2a} this means that
$J(x)$ has rank $m-1$ for every $0 \neq x \in X_{m,n}$. Thus every 
$0 \neq x \in X_{m,n}$ is smooth.
\end{pf}

\begin{lem} \label{lem7.2c}
If $1 \leq m \leq n-1$ 
then $X_{m, n}$ is an irreducible variety of dimension $n-m+1$.
\end{lem}

\begin{pf}
Consider the morphism $\pi \colon X_{m,n} \lra \bbA^{n-m+1}$ given by
\begin{equation} \label{e7.15}
\pi(x) = (s_m(x), \ldots, s_n(x)) \; 
\end{equation} 
where $s_j$ is the $j$th elementary symmetric polynomial, as before.
Then $\pi$ is
surjective, and the fibers of $\pi$ are precisely the $\Sym_n$-orbits
in $X_{m,n}$. This shows that $\dim(X_{m,n}) = n-m+1$.  On the other hand,
since $X_{m,n}$ is cut out by $m-1$ homogeneous polynomials in $\bbA^n$, 
every irreducible component of it has dimension $\geq n-m+1$; cf., 
e.g.,~\cite[Proposition I.7.1]{Hart}. We conclude that every component
of $X_{m, n}$ has dimension exactly $n-m+1$ and the restriction of
$\pi$ to any component of $X_{m,n}$ is dominant. Since
$\Sym_n$ acts transitively on the fibers of $\pi$, its action on the set of 
the irreducible components of $X_{m,n}$ is also transitive. 

Let $X_1$ be an irreducible component of $X_{m,n}$ and let $H$ be the subgroup
of $\Sym_n$ preserving $X_1$. Since $\Sym_n$ transitively permutes the
components of $X_{m,n}$, it is enough to show that $H = \Sym_n$. We will
do this by proving that $H$ contains every transposition $(i, j)$ for
$1 \leq i < j \neq n$.

We claim that $\Stab(x) \subset H$ for every $0 \neq x \in X_1$. Indeed,
assume to the contrary that $g \in \Stab(x)$ but $g(X_1) \neq X_1$. 
Then $g(X_1)$ and $X_1$ are distinct irreducible components of $X_{m,n}$ 
passing through $x$. Hence, $x$ is a singular point of $X_{m,n}$, 
contradicting Lemma~\ref{lem7.2b}. This proves the claim.

It is now sufficient to show that for every transposition 
$g = (i, j)$ there exists a point $0 \neq x \in X_1$ such that $g(x) = x$. 
In other words, we want to show
that there is a non-zero point $x=(x_1, \ldots, x_n) \in X_1$ with $x_i = x_j$. 

To prove the last assertion, we pass to the projective space $\bbP^{n-1}$. Let
$\bbP(X_{m,n})$ be the projectivization of $X_{m,n}$, i.e., 
the subvariety of $\bbP^{n-1}$ given by \eqref{e7.1}. Then the irreducible
components of $X_{m,n}$ are affine cones over the irreducible components
of $\bbP(X_{m,n})$; in particular, $X_1$ is an affine cone over $\bbP(X_1)$,
where $\dim(\bbP(X_1))=\dim(X_1) - 1 = n-m$. Thus our assumption that
$m \leq n-1$ translates into $\dim(\bbP(X_1)) \geq 1$. Thus $\bbP(X_1)$ has
a non-trivial intersection with any hyperplane. In particular,
$\bbP(X_1) \cap \{ x_i = x_j \} \neq \emptyset$ and, hence,
$X_1$ contains a non-zero point preserved by $(i, j)$.
This completes the proof of Lemma~\ref{lem7.2c}.
\end{pf}

\begin{remark} \label{rem.x_n,n} 
The condition $m \leq n-1$ in Lemma~\ref{lem7.2c} is essential. Indeed, 
the variety $X_{n,n}$ is a union of $(n-1)!$ lines
given (in parametric form) by $(\zeta_1t, \zeta_2t, \dots, \zeta_nt)$, where
$\zeta_1, \dots, \zeta_n$ are distinct $n$-th roots of unity. In other words,
$\bbP(X_{m,n})$ is a union of the $(n-1)!$ projective points of the form
$(\zeta_1: \dots :\zeta_n)$; note that none of these
points lies on the hyperplane $x_i = x_j$ for any choice of 
$1 \leq i < j \leq n$.
\end{remark}

\subsection*{Proof of Theorem~\ref{thm7.1}}

To prove the inequality $\tau(E/F) \leq n-m$, let 
$z = \dfrac{a_{n-1}}{a_n} x$.
(Note that here we are using the assumption
$m \leq n-1$.) Substituting 
$x = \dfrac{a_n\mathstrut}{a_{n-1\mathstrut}} z$ into the equation
$f(x) = 0$, we see that
the minimal polynomial of $z$ over $F$ is of the form
$$f_z(t) = t^n + b_m t^m + \dots + b_{n-1} t + b_n \; , $$
where $b_n = b_{n-1} = \dfrac{a_{n-1}^n}{a_n^{n-1}}$. Thus
$$ \tau(E/F) \leq \trdeg_k k(b_m, \dots, b_{n-1}, b_n) =
\trdeg_k k(b_m, \dots, b_{n-1}) \leq n-m \; ,$$
as claimed.

It therefore remains to show that $\tau(E/F) \geq n-m$.
Since \[ \tau(E \otimes_k \overline{k}/F \otimes_k \overline{k}) 
\geq \tau(E/F)\, , \] we may assume 
without loss of generality that $k = \overline{k}$ is algebraically 
closed; cf.\ Remark~\ref{rem6.05}.
Let  $ X_{m, n}$ be the $\Sym_n$-variety defined by \eqref{e7.1} and
let $E^{\#}$ be the normal closure of $E$ over $F$. Note that by 
\cite[Lemma 2.3]{br1} $\tau(E/F) = \tau(E^{\#}/F)$. Our strategy will
thus be as follows: first we will show that 
\begin{equation} \label{e7.2}
\tau(E^{\#}/F) = \ed(X_{m, n}) \; , 
\end{equation}
then
\begin{equation} \label{e7.3}
\ed(X_{m, n}) \geq n-m \; .
\end{equation}
We now proceed to prove \eqref{e7.2}. By \cite[Lemma 2.7]{br1} it is 
enough to show that the field extensions 
$E^{\#}/F$ and $k(X_{m, n})/k(X_{m, n})^{\Sym_n}$ are isomorphic. 

We claim that $k(X_{m, n})^{S_n} = k(s_m, \ldots, s_n)$, where $s_i$ is the
$i$th symmetric polynomial of $x_1, \ldots, x_n$, viewed as a regular
function on $X_{m, n}$. Indeed, it is clear that
$k(s_1, \ldots, s_n) \subset k(X_{m, n})^{S_n}$. To prove equality,
observe that the polynomial 
$f(x) = x^n + a_m x^{n-m} + \dots + a_{n-1} x + a_n$ 
has $n$ distinct roots for a generic choice of 
$(a_m, a_{m+1}, \dots, a_n) \in \bbA^{n-m+1}$ (because $x^n -1$ has
$n$ distinct roots).  This means that the map
\[
\pi\colon X_{m, n} \lra \bbA^{n-m+1}
\]
given by \eqref{e7.15}, is generically $n!:1$ and consequently,
$[k(X_{m, n}): k(s_m, \ldots, s_n)] = n!$. Thus
\[ [k(X_{m, n})^{S_n} : k(s_m, \ldots, s_n)] = \dfrac
{[k(X_{m, n}) : k(s_m, \ldots, s_n)]} {[k(X_{m, n}) : k(X_{m, n})^{S_n}]}
= \dfrac{n!}{n!} = 1 \; , \]
as claimed.

Continuing with the proof of \eqref{e7.2},
note that the $s_m, \ldots, s_n$ are algebraically independent
over $k$. (This follows, e.g., from the fact that
the map $\pi$ defined in \eqref{e7.15}, is dominant.) Thus the 
fields $k(s_m, \ldots, s_n)$ and $F = k(a_m, \ldots, a_n)$ are isomorphic
via a map that takes $s_i$ to $a_i$ for every $i$. Now observe that
$k(X_{m, n})$ is the splitting field of the polynomial
$g(x) = x^n + s_m x^{n-m} + \ldots + s_{n-1}x + s_n$ over $k(X_{m, n})^{\Sym_n}
=k(s_m, \ldots, s_n)$ and $E^{\#}$ is by definition the splitting field of
$f(x)$ over $F = k(a_m, \ldots, a_n)$. By the uniqueness of the splitting
field, we see that the field extensions $k(X_{m, n})/k(X_{m, n})^{\Sym_n}$
and $E^{\#}/F$ are isomorphic, as claimed.
This completes the proof of~\eqref{e7.2}.

It remains to prove the inequality \eqref{e7.3}. 
In view of Theorem~\ref{thm6.1}(2) it is sufficient to show that
there exists a smooth point $x \in X_{m, n}$ such that $\Stab(x)$ contains 
a subgroup isomorphic to $(\bbZ/2\bbZ)^{n-m}$. 
We shall thus look for a point of the form
\begin{equation} \label{e7.4}
x = (\alpha_1, \alpha_1, \alpha_2, \alpha_2, \dots, \alpha_{n-m}, 
\alpha_{n-m}, \alpha_{n-m+1}, \alpha_{n-m+2}, \ldots, \alpha_{m-1}, \alpha_m)
\; , 
\end{equation}
where at least one $\alpha_i$ is non-zero. (Here we are using the assumption
that $m \geq n/2$ and thus $2(n-m) \leq n$.) 
By Lemma~\ref{lem7.2b} any non-zero point $x$ of $X_{m, n}$ is smooth;
moreover, if $x$ is as in \eqref{e7.4} then
$\Stab(x)$ contains the subgroup 
\[ \lf< (1, 2) \, , \; (3, 4) \, , \; \dots \, , \; (2n-2m-1, 2n-2m) \r> 
\simeq (\bbZ/2\bbZ)^{n-m} \; . \] 

Thus we only need to show that a non-zero point of the form \eqref{e7.4}
exists on $X_{m, n}$.  Substituting $x$ into the defining equations
$p_1(x) = \ldots = p_{m-1}(x) = 0$ of $X_{m, n}$ (see \eqref{e7.1-1}),
we obtain 
a system of $m-1$ homogeneous equations in $\alpha_1, \ldots, \alpha_m$.
Since the number of variables is greater than the number of equations,
this system has a non-trivial solution, which gives us 
the desired point. This completes the proof of the inequality \eqref{e7.3}
and, hence, of Theorem~\ref{thm7.1}.
\qed

\begin{remark} \label{rem7.3} 
The same argument (with part (3) of Theorem~\ref{thm6.1} used in place of
part (2)) shows that $\tau(E; 2) = n-m$ in the sense of 
Definition~\ref{def.tau} (here, as before, $E$ is viewed as an $n$-dimensional
$F$-algebra).  In particular, the polynomial $f(x)$ of Theorem~\ref{thm7.1}
cannot be reduced to a form with $\leq n-m$ algebraically independent 
coefficients by a Tschirnhaus transformation, even if we allow auxiliary
radicals of odd degree; cf.~\cite[Theorem 7.1]{br2}.
\end{remark}

\smallskip

\noindent Department of Mathematics, Oregon State University,
Corvallis, OR 97331-4506, USA

\noindent Current mailing address: 
PMB 136, 333 South State St., Lake Oswego, OR 97034-3961, USA.

\noindent email: zinovy@@math.orst.edu 

\smallskip

\noindent
Department of Mathematics and Computer Science,
University of the Negev, Be'er Sheva', Israel

\noindent
Current mailing address: Hashofar 26/3, Ma'ale Adumim, Israel.

\noindent
email: youssin@@math.bgu.ac.il

\bigskip
\bigskip
%%%%%%%%%%%%%%%%%%%%%%%%%%%

% Here begins the appendix, by Koll\'ar and Szab\'o.

\begingroup % make the definitions local for the appendix only

% Some versions of amslatex have \mathbb and some do not.  The following code
% defines \mathbb to be \Bbb if the former is undefined.

\ifx\mathbb\undefined \newcommand{\mathbb}{\Bbb}\else\fi

% If this does not work, the following line, if uncommented, should give \mathbb
% the value of \Bbb in any case:
% \let\mathbb=\Bbb

\let \cedilla =\c
\renewcommand{\c}[0]{{\mathbb C}}  
\let \crossedo =\o
\renewcommand{\o}[0]{{\cal O}} 
\newcommand{\z}[0]{{\mathbb Z}}
\newcommand{\n}[0]{{\mathbb N}}

% \let \ringaccent=\r  %%% \r shorthand for 'ring accent' % this line has been
% moved to the preamble
%\renewcommand{\r}[0]{{\mathbb R}} % our version of amslatex does not have \r
\renewcommand{\r}[0]{{\mathbb R}} %\r has different meaning in the other
                                  % parts of the paper

\renewcommand{\a}[0]{{\mathbb A}}

\newcommand{\p}[0]{{\mathbb P}}
\newcommand{\f}[0]{{\mathbb F}}
\newcommand{\q}[0]{{\mathbb Q}}
\newcommand{\map}[0]{\dasharrow}
\newcommand{\qtq}[1]{\quad\mbox{#1}\quad}
\newcommand{\spec}[0]{\operatorname{Spec}}
\newcommand{\pic}[0]{\operatorname{Pic}}
\newcommand{\gal}[0]{\operatorname{Gal}}
\newcommand{\cont}[0]{\operatorname{cont}}
\newcommand{\mult}[0]{\operatorname{mult}}
\newcommand{\discrep}[0]{\operatorname{discrep}}
\newcommand{\totaldiscrep}[0]{\operatorname{totaldiscrep}}

\newcommand{\supp}[0]{\operatorname{Supp}}    
\newcommand{\red}[0]{\operatorname{red}}    
\newcommand{\im}[0]{\operatorname{im}}    
\newcommand{\proj}[0]{\operatorname{Proj}}    
\newcommand{\wt}[0]{\operatorname{wt}}    
\newcommand{\socle}[0]{\operatorname{socle}}    
\newcommand{\coker}[0]{\operatorname{coker}}    
\newcommand{\ext}[0]{\operatorname{Ext}}    
\renewcommand{\trace}[0]{\operatorname{Trace}}   % already in the preamble
\newcommand{\cent}[0]{\operatorname{center}}
\newcommand{\bs}[0]{\operatorname{Bs}}
\newcommand{\specan}[0]{\operatorname{Specan}}    
\newcommand{\inter}[0]{\operatorname{Int}}    
\newcommand{\sing}[0]{\operatorname{Sing}}    
\newcommand{\ex}[0]{\operatorname{Ex}}    
\newcommand{\chr}[0]{\operatorname{char}}    

\newcommand{\nec}[1]{\overline{NE}({#1})}

\newcommand{\rup}[1]{\lceil{#1}\rceil}
\newcommand{\rdown}[1]{\lfloor{#1}\rfloor}

\def\into{\DOTSB\lhook\joinrel\rightarrow}

\renewcommand{\thesection}{A}

\section{Appendix}
\vbox{
\begin{center} {\LARGE {\bf Fixed Points of Group Actions and Rational Maps}} 
\end{center}
\smallskip

\begin{center} by \end{center}

\begin{center}
{\Large J{\'a}nos Koll{\'a}r}\\ 
(Department of Mathematics, Princeton University, Princeton, \\
NJ 08544-1000, USA, kollar@@math.utah.edu)
\end{center}

\begin{center} and \end{center}

\begin{center}
{\Large Endre Szab{\'o}}\\ (Mathematical Institut, 
Budapest, PO.Box 127, 1364 Hungary, endre@@math-inst.hu)
\end{center}
}

\bigskip
The aim of this note is to give simple proofs of
the results in Section~\ref{sect5}
about the behaviour of fixed points of finite group
actions under rational maps. Our proofs work in any
characteristic.

\begin{lem}\label{class.defn}
 Let $K$ be an algebraically closed field and  $H$
a (not necessarily connected) linear
algebraic group over $K$.  The following are equivalent.
\begin{enumerate}
\item Every representation  $H\to GL(n,K)$ has an
$H$-eigenvector.
\item There is a  (not necessarily connected) unipotent, normal 
subgroup $U<H$ such that $H/U$ is abelian.
\end{enumerate}
\end{lem}

Proof. Let $H\to GL(n,K)$ be a faithful representation.
If (\ref{class.defn}.1) holds then $H$ is conjugate to
an upper triangular subgroup, this implies (\ref{class.defn}.2).

Conversely, any representation of a unipotent group
has fixed vectors (cf.\ \cite[I.4.8]{borel-book})
and the subspace of all fixed vectors is an
$H/U$-representation.\qed

\begin{prop}[Going down]\label{Going-down}
  Let $K$ be an algebraically closed
field,
$H$ a linear algebraic group over $K$ and $f:X\map Y$ an
$H$-equivariant map of
$K$-schemes. Assume that
\begin{enumerate}
\item $H$ satisfies the equivalent conditions
of (\ref{class.defn}),
\item $H$ has a smooth fixed point on $X$, and
\item $Y$ is proper.
\end{enumerate}
\noindent Then $H$ has a fixed point on $Y$.
\end{prop}

Proof. The proof is by induction on $\dim X$. The  case  $\dim
X=0$ is clear.

Let $x\in X$ be a smooth $H$-fixed point and consider the blow up
$B_xX$ with exceptional divisor $E\cong \p^{n-1}$. 
The $H$-action lifts to $B_xX$ and so we get an $H$-action on $E$
which has a fixed point by (\ref{class.defn}.1).  Since $Y$ is
proper, the induced rational map $B_xX\to X\map Y$ is defined
outside a subset of codimension at least 2. Thus we get an
$H$-equivariant rational map $E\map Y$. By induction, there is a 
fixed point on $Y$.\qed

\begin{remark} \label{rem.Going-down} 
If $H$ does not satisfy the  conditions
of (\ref{class.defn})
then (\ref{Going-down}) fails for some actions.
Indeed, let $H\to GL(n,K)$ be a representation without an 
$H$-eigenvector.  This gives an $H$-action on $\p^n$ with a single fixed
point $Q\in \p^n$. The corresponding action on $B_Q\p^n$ has no fixed
points.
\end{remark}

\begin{prop}[Going up]\label{Going-up}
  Let $K$ be an algebraically closed field
and
$H$ a finite abelian group of prime power order $q^n$ 
($q$ is allowed to coincide with $\chr K$). Let $p:X\map Z$ be an
$H$-equivariant map of irreducible
$K$-schemes. Assume that
\begin{enumerate}
\item $p$ is generically finite, dominant and $q\not\vert\deg (X/Z)$,
\item $H$ has a smooth fixed point on $Z$, and
\item $X$ is proper.
\end{enumerate}
\noindent Then $H$ has a fixed point on $X$.
Moreover, if $X\map Y$ is an $H$-equivariant
map to a proper $K$-scheme then $H$ has a fixed point on $Y$.
\end{prop}

Proof.  The proof is by induction on $\dim Z$. 
The  case  $\dim Z=0$ is clear.

Let $z\in Z$ be a smooth fixed point and $E\subset B_zZ$
the exceptional divisor. Let $\bar p:\bar X\to B_zZ$
denote the normalization of $B_zZ$ in the field of rational
functions of $X$ and $F_i\subset \bar X$ the divisors lying over
$E$.  $H$ acts on the set   $\{F_i\}$. Let ${\cal F}_j$
denote the $H$-orbits and in each pick a divisor
$F^*_j\in {\cal F}_j$.
By the ramification formula (see \cite[Corollary XII.6.2]{lang65})
$$
\deg (X/Z)=\sum_j 
|{\cal F}_j|\cdot\deg (F^*_j/E)\cdot e(\bar p, F^*_j)
$$
where $e(\bar p,F^*_j)$ denotes
the ramification index of $\bar p$ at the generic point of $F^*_j$.
Since $\deg (X/Z)$ is not divisible by $q$,
there is an orbit ${\cal F}_0$ consisting of  a single element
$F^*_0$  such that $\deg (F^*_0/E)$ is not divisible by $q$.

We have $H$-equivariant rational maps 
$F^*_0\map E$, $F^*_0\map X$ and $F^*_0\map Y$.
By induction $H$ has a fixed point on $F_0^*$, $X$ and $Y$.
\qed

\begin{remark} We see from the proof that (\ref{Going-up})
also holds  if $H$ is abelian and only one of the prime divisors of
$|H|$ is less than $\deg (X/Z)$.
\end{remark}

The method  also gives a simpler proof of a result
of \cite{nishi}. One can view this as a version of
(\ref{Going-down}) where $H$ is the absolute Galois group of $K$.

\begin{prop}[Nishimura lemma]\label{nishimura}
  Let $K$ be a field  and $f:X\map Y$ a rational map of
$K$-schemes. Assume that
\begin{enumerate}
\item $X$ has a smooth $K$-point, and
\item $Y$ is proper.
\end{enumerate}
\noindent Then $Y$ has a   $K$-point.
\end{prop}

Proof. The proof is by induction on $\dim X$. The  case  $\dim
X=0$ is clear.

Let $x\in X$ be a smooth $K$-point and consider the blow up
$B_xX$ with exceptional divisor $E\cong \p^{n-1}$. The divisor 
$E$ has smooth $K$-points.  Since $Y$ is
proper, the induced rational map $B_xX\to X\map Y$ is defined
outside a subset of codimension at least 2 and we get a rational
map $E\map Y$. By induction, there is a  $K$-point on $Y$.\qed

\begin{remark}  One can combine (\ref{Going-down}) and 
(\ref{nishimura}) if we know that any $H$-representation has an
eigenvector defined over $K$. There are two interesting cases
where this condition holds:
\begin{enumerate}
\item  $H$ is Abelian of order $n$ and $K$ contains all $n$th
roots of unity.
\item $H$ is nilpotent and its order is a power of $\chr K$.
\end{enumerate}
\end{remark}

\endgroup  % end of the appendix
\end{document}